%
%
%
%
%
%
%
\input amstex
\documentstyle{amsppt}
\NoBlackBoxes

\topmatter
\title Recent Developments in Toric Geometry\endtitle
\author David A. Cox\endauthor
\leftheadtext{DAVID A. COX}%

\address Department of Mathematics and Computer Science, Amherst
College, Amherst, Massachusetts 01002 \endaddress

\email dac\@cs.amherst.edu\endemail


\subjclass Primary 14M25\endsubjclass

\abstract This paper will survey some recent developments in
the theory of toric varieties, including new constructions of toric
varieties and relations to symplectic geometry, combinatorics and
mirror symmetry.\endabstract

\thanks The author was supported in part by NSF Grant
\#DMS--9301161.\endthanks 

\endtopmatter

\document

\head Introduction\endhead

	A toric variety over ${\Bbb C}$ is a $n$-dimensional normal
variety $X$ containing $({\Bbb C}^*)^n$ as a Zariski open set in such
a way that the natural action of $({\Bbb C}^*)^n$ on itself extends to
an action of $({\Bbb C}^*)^n$ on $X$.  This seemingly simple
definition leads to a fascinating combinatorial structure and some
surprisingly rich mathematics.  In this article, we will discuss some
recent developments in toric varieties, including novel applications
and new foundations for the entire theory.

	Toric varieties were discovered in the early 1970's
independently by several groups of people.  From the beginning, the
theory of toric varieties has led to some notable applications,
including the following:
\roster
\item"{$\bullet$}" The characterization of algebraic subgroups of
maximal rank of Cremona groups, in Demazure's 1970 paper \cite{Dem}.
\item"{$\bullet$}" The stable reduction theorem, in the 1973 book
\cite{KKMS} by Knudsen, Kempf, Mumford and Saint-Donat.
\item"{$\bullet$}" The construction of nice (meaning ``toroidal'')
compactifications of discrete quotients of bounded symmetric domains,
begun in the 1973 paper \cite{Sat} by Satake and the 1975 book
\cite{AMRT} by Ash, Mumford, Rapaport and Tai.
\item"{$\bullet$}" The rich connections between Newton polytopes,
toric varieties and singularities, first explored by Kushnirenko
\cite{Kus} in 1976 and Khovanskii \cite{Kho} in 1977.
\item"{$\bullet$}" The use of Hard Lefschetz for simplicial toric
varieties to prove McMullen's conjectures for the number of vertices, 
edges, faces, etc.~of convex simplical polytopes, in Stanley's 1980
paper \cite{Sta}.
\endroster
\noindent This brief list, of course, does not do justice to the work
of many other people who have written about toric varieties.

	Besides these applications of toric varieties, many people
have come to the realization, as noted by Fulton in his 1993 book
\cite{Ful}, that ``toric varieties have provided a remarkably fertile
testing ground for general theories.''  An example of this philosophy
can be found in Reid's 1983 paper \cite{Rei2} which studies Mori theory
in the context of toric varieties.

	The theory of toric varieties has also had the benefit of some
superb exposition.  We've already mentioned Fulton's book \cite{Ful},
and two other classic references in the field are Danilov's 1978
survey \cite{Dan2} and Oda's 1988 book \cite{Oda2}.  We warmly recommend
these to anyone who wants to learn more about this fascinating part of
algebraic geometry.

	Although the classical theory seems fairly complete, the last
few years have seen an explosion of new ideas and tools for studying
toric varieties as well as significant new applications and relations
to other areas of mathematics (and physics!).  This paper will survey
some of these new developments.  We begin in \S1 with a review of the
notation and terminology we will use.  Then \S2 introduces a new
construction of toric varieties, similar to way that ${\Bbb P}^n$ is
realized as a quotient of ${\Bbb C}^{n+1}-\{0\}$.  It follows that
simplicial toric varieties have ``homogeneous coordinates'' which
enable us to define subvarieties by global equations.  In \S3, we
briefly discuss the K\"ahler cone of a toric variety.  Then \S4
describes a related construction of smooth simplicial toric varieties
coming from symplectic geometry.

	In \S5, we discuss a yet another method for constructing toric
varieties which has gained prominence recently.  Here, a toric variety
is defined to be the closure of an equivariant map from a torus.  Such
toric varieties need not be normal, which gives the theory a slightly
different flavor.  Then \S6 shows how classical results of Griffiths
on the cohomology of projective hypersurfaces can be generalized to
the toric case.  The {\it secondary fan\/} is the subject of \S7.
This fan, which has some remarkable applications, begins with a
collection of rays and asks how many fans have these rays as their
$1$-dimensional cones.  In \S8 we will discuss reflexive polytopes and
Calabi-Yau hypersurfaces, and then \S9 will touch briefly on the work
of Gelfand, Kapranov and Zelevinsky on resultants, discriminants and
hypergeometric equations.

	Mirror symmetry is the topic of \S10.  This field represents a
fascinating interaction between mathematics and physics.  Toric
geometry plays a prominent role in many mirror symmetry contructions,
and there are even physical theories specially designed for toric
varieties.  In particular, we will see that mirror symmetry makes use
of virtually {\it everything\/} in \S\S2--9.  Finally, \S11 will
discuss some other research being done on toric varieties.

	In this survey, we will assume that the reader is familiar
with basic theory of toric varieties as presented in \cite{Dan2},
\cite{Ful} or \cite{Oda2}.  For simplicity, we will usually work over
the complex numbers ${\Bbb C}$.  

	To keep the bibliography from getting too large, references to
some topics are not complete---it was often more convenient to refer
to later papers rather than the original ones.  A fuller picture of
recent work on toric varieties can be obtained by checking the
references in the papers mentioned in this survey.  The reader may
also want to consult the 1989 toric survey of Oda \cite{Oda3}, which
has an extensive bibliography.

\head \S1. Notation and terminology for toric varieties\endhead

	The basic combinatorial object associated to a toric variety
is a fan.  One starts with an integer lattice $M \simeq {\Bbb Z}^n$
and its dual lattice $N$.  Then a {\it fan\/} $\Sigma$ in $N_{{\Bbb
R}} = N\otimes {\Bbb R}$ consists of a finite collection of {\it
strongly convex rational polyhedral cones\/} $\sigma \subset N_{{\Bbb
R}}$ which is closed under intersection and taking faces.  The
$1$-dimensional cones $\rho$ play a prominent role in the theory, and
it is customary to denote the unique generator of $\rho\cap N$ by the
same letter $\rho$.  Finally, we say that $\Sigma$ is {\it
simplicial\/} if the minimal generators of each $\sigma \in \Sigma$
are linearly independent over ${\Bbb R}$.  (For details of these and
other definitions in this section, consult the references given in the
introduction.)

	\subhead The classical construction\endsubhead Given a fan
$\Sigma$, each cone $\sigma \in \Sigma$ has a {\it dual
cone\/} $\sigma^\vee = \{v \in M_{{\Bbb R}}: \langle v,\sigma\rangle
\ge 0\}$, which determines the semigroup algebra
${\Bbb C}[\sigma^\vee\cap M]$.  In the classical formulation of the
theory, the toric variety $X = X_\Sigma$ is obtained from $\Sigma$ by
gluing together the affine toric varieties $X_\sigma =
\roman{Spec}({\Bbb C}[\sigma^\vee\cap M])$ for $\sigma \in \Sigma$.
We will see in~\S2 that there are now other ways to construct $X$ from
$\Sigma$.

	The torus $T = N\otimes{\Bbb C}^*$ sits inside $X$ as
described in the introduction.  Hence $N$ gives the $1$-parameter
subgroups of $T$ and $M$ is the character group.  So each $m \in M$
gives $\chi^m : T \to {{\Bbb C}}^*$, which can be regarded as a
rational function on $X$.

	Although a toric variety $X$ can be singular, it is always
Cohen-Macaulay (this is proved in \cite{Hoc} and \cite{Dan2}), so that
duality theory works nicely.  In particular, the dualizing sheaf
$\omega_X$ coincides with the sheaf $\widehat\Omega^n_X$ of Zariski
$n$-forms on $X$.  An especially nice case is when $X$ is simplicial
(meaning that $\Sigma$ is simplicial).  In this case, $X$ is a
$V$-manifold and for many purposes (including cohomology over ${\Bbb
Q}$ and Hodge theory) behaves like a manifold.

	\subhead The role of polyhedra\endsubhead One way to see the
connection with polyhedra is via divisors on $X$.  Each
$1$-dimensional cone $\rho \in \Sigma$ corresponds to a Weil divisor
$D_\rho \subset X$.  A divisor $D = \sum_\rho a_\rho D_\rho$ gives a
(possibly unbounded) convex polyhedron
$$
\Delta_D = \{m \in M_{{\Bbb R}}: \langle m,\rho\rangle \ge -a_\rho\}
\subset M_{{\Bbb R}}.\tag1.1
$$
To see how $\Delta_D$ relates to $D$, consider the reflexive sheaf
$\Cal{O}_X(D)$ whose sections over $U \subset X$ are those rational
functions $f$ such that $\roman{div}(f) + D \ge 0$ on $U$.  Then
$$
H^0(X,\Cal{O}_X(D)) = \bigoplus_{m \in \Delta_D\cap M} {\Bbb C}\cdot
\chi^m. \tag1.2
$$
When $D = \sum_\rho a_\rho D_\rho$ is Cartier, there is a {\it support
function\/} $\psi_D$ with the property that for each $\sigma \in
\Sigma$, there is $m_\sigma \in M$ such that $\psi_D(\rho) = \langle
m_\sigma,\rho\rangle = -a_\rho$ for all $\rho \subset \sigma$.  In
particular, $\psi_D$ is linear and integral on each cone in $\Sigma$.

	When $X$ is complete, a classic fact is that $D$ is ample if
and only if $\psi_D$ is strictly convex.  In this case, $\Delta_D$ is
a $n$-dimensional integral convex polytope (= bounded polyhedron)
which is combinatorially dual to $\Sigma$, i.e., facets of $\Delta_D$
(faces of dimension $n-1$) correspond to $1$-dimensional cones $\rho
\in \Sigma$ and, more generally, $i$-dimensional faces of $\Delta_D$
correspond to $(n-i)$-dimensional cones of $\Sigma$.  

	Conversely, given a $n$-dimensional integral convex polytope
$\Delta \subset M_{{\Bbb R}}$, there is a unique fan (sometimes called
the {\it normal fan\/} of $\Delta$) such that corresponding toric
variety $X_\Delta$ has a Cartier divisor which gives $\Delta$ exactly
(see \cite{Oda1{\rm, Sect.~2.4}}).  In \S5 we will give a method (due
to Batyrev) for obtaining $X_\Delta$ directly from $\Delta$.

	\subhead Other remarks\endsubhead A minor omission in the
classic references for toric varieties is that line bundles are
discussed in detail, but not the reflexive sheaves coming from Weil
divisors $\sum_\rho a_\rho D_\rho$.  Basic facts about reflexive
sheaves on normal varieties can be found in \cite{Rei1}.  We should
also mention that there is a standard conflict of notation: some
authors use $\Delta$ for the fan (see \cite{Ful},
\cite{Oda2}, \cite{Stu2}), while others use $\Sigma$ for the fan and
$\Delta$ for a polytope (see \cite{Dan2}, \cite{Bat1}).

	It is also possible to consider infinite fans.  Gluing
together affine toric varieties for cones in such a fan would lead to
a scheme which is only locally of finite type.  However, given an
appropriate discrete group action, one can then take a quotient to get
a variety.  This is the approach used in \cite{KMMS} and \cite{AMRT}.
A nice example is the resolution of a $2$-dimensional Hilbert cusp
singularity, which is discussed in \cite{Oda1{\rm, Sect.~4.1}}.

\head \S2. Global coordinates for toric varieties\endhead

	Projective space is one of the simplest examples of a toric
variety.  The way ${\Bbb P}^n$ is obtained by gluing together affine
spaces is a special case of the classic construction of a toric
variety.  But ${\Bbb P}^n$ can also be constructed as a quotient
$$
{\Bbb P}^n = ({\Bbb C}^{n+1} - \{0\}) / {\Bbb C}^*,
$$
which is where we get the homogeneous coordinates on projective
space.  Recently, this construction has been generalized to most toric
varieties.  We will describe this construction and some of its
consequences. 

	\subhead The construction\endsubhead We begin by fixing a fan
$\Sigma$ in $N_{{\Bbb R}} \simeq {\Bbb R}^n$ and letting $\Sigma(1)$
be the set of $1$-dimensional cones in $\Sigma$.  We will assume that
$$
\hbox{the $1$-dimensional cones $\rho \in \Sigma(1)$ span $N_{{\Bbb
R}}$} \tag2.1
$$
(where as usual we regard $\rho$ as the integral generator of its
cone).  Any complete fan satisfies this condition.  Then consider the
affine space ${\Bbb C}^{\Sigma(1)}$ with variables $x_\rho$ for $\rho
\in \Sigma(1)$.  We need to remove a certain {\it exceptional
subset\/} from ${\Bbb C}^{\Sigma(1)}$.  This is defined as follows.
For each $\sigma \in \Sigma$, consider the monomial $\hat x_\sigma =
\Pi_{\rho \not\subset\sigma} x_\rho$.  These generate the monomial
ideal $B(\Sigma) = \langle \hat x_\sigma : \sigma \in
\Sigma\rangle$, and the exceptional set $Z(\Sigma) \subset {\Bbb
C}^{\Sigma(1)}$ is the subvariety defined by $B(\Sigma)$.

	The toric variety $X = X_\Sigma$ will be a quotient of ${\Bbb
C}^{\Sigma(1)} - Z(\Sigma)$ by the group 
$$
G = \roman{Hom}_{{\Bbb Z}}(A_{n-1}(X),{\Bbb C}^*),
$$
where $A_{n-1}(X)$ is the Chow group of Weil divisors modulo rational
equivalence.  To see how this group acts on ${\Bbb C}^{\Sigma(1)} -
Z(\Sigma)$, recall the exact sequence
$$
0 \longrightarrow M \overset\alpha\to\longrightarrow \bigoplus_{\rho
\in \Sigma(1)}{\Bbb Z}\cdot D_\rho \overset\beta\to\longrightarrow
A_{n-1}(X) \longrightarrow 0 \tag2.2
$$
where $\alpha$ sends $m \in M$ to $\roman{div}(\chi^m) = \sum_\rho
\langle m,\rho\rangle D_\rho$ and $\beta$ is the obvious map from Weil
divisors to the Chow group.  The injectivity of $\alpha$ is equivalent
to \thetag{2.1}.  Applying $\roman{Hom}_{{\Bbb Z}}(-,{\Bbb C}^*)$
yields the exact sequence
$$
1 \longrightarrow G \longrightarrow ({\Bbb C}^*)^{\Sigma(1)}
\longrightarrow N\otimes{\Bbb C}^* \longrightarrow 1.\tag2.3 
$$
This gives a natural action of $G$ on ${\Bbb C}^{\Sigma(1)}$, and
since $Z(\Sigma)$ is a union of coordinate subspaces, $G$ preserves
${\Bbb C}^{\Sigma(1)} -Z(\Sigma)$.

	A nice example of what this looks like is given by ${\Bbb
P}^n$.  The fan $\Sigma$ of ${\Bbb P}^n$ is well-known, and we leave
it to the reader to check that in this case, the monomial ideal
$B(\Sigma)$ is the ``irrelevant'' ideal $\langle
x_0,\ldots,x_n\rangle$, so that the exceptional set $Z(\Sigma)$
consists of the origin.  Furthermore, the exact sequence \thetag{2.2}
becomes
$$
0 \longrightarrow {\Bbb Z}^n \longrightarrow {\Bbb Z}^{n+1}
\overset\beta\to\longrightarrow {\Bbb Z} \longrightarrow 0,
$$
where $\beta(a_0,\ldots,a_n) = \sum_{i=0}^n a_i$.  By
\thetag{2.3}, the action of $G$ on ${\Bbb C}^{\Sigma(1)}
-Z(\Sigma)$ is the usual action of ${\Bbb C}^*$ on ${\Bbb C}^{n+1} -
\{0\}$.  The quotient, of course, is ${\Bbb P}^n$.

	Returning to the general case, assume for the moment that $X$
is given by the quotient $\big({\Bbb C}^{\Sigma(1)}
-Z(\Sigma)\big)/G$.  Then the natural inclusion $ ({\Bbb
C}^*)^{\Sigma(1)} \subset {\Bbb C}^{\Sigma(1)} -Z(\Sigma)$, combined
with \thetag{2.3}, gives an inclusion $T = N\otimes{\Bbb C}^* \subset
X_\Sigma$ as a dense open set.  The action of $T$ on $X$ is also easy
to see in this picture: it is inherited from the action of the ``big
torus'' $({\Bbb C}^*)^{\Sigma(1)}$ on ${\Bbb C}^{\Sigma(1)}
-Z(\Sigma)$.

	Taking quotients in algebraic geometry can be a bit subtle,
and the ``quotient'' $\big({\Bbb C}^{\Sigma(1)} -Z(\Sigma)\big)/G$ is
no exception.  The precise relation between this quotient and the
toric variety is as follows.

\proclaim{Theorem 2.1} Let $X$ be a toric variety whose fan $\Sigma$
satisfies \thetag{2.1}.  Then:
\roster
\item $X$ is the universal categorical quotient $\big({\Bbb
C}^{\Sigma(1)} -Z(\Sigma)\big)/G$.  
\item $X$ is a geometric quotient  $\big({\Bbb C}^{\Sigma(1)}
-Z(\Sigma)\big)/G$ if and only if $\Sigma$ is simplicial. 
\endroster
\endproclaim

\demo{Remarks on Theorem 2.1} A universal categorical quotient is
a $G$-equi\-variant map $\pi : {\Bbb C}^{\Sigma(1)} -Z(\Sigma) \to X$
(where the action on $X$ is trivial) which is universal in the obvious
sense.  It is a geometric quotient precisely when the fibers of $\pi$
coincide with the $G$-orbits.  Hence the simplicial case is the one
closest to the way we think about projective space.

	To get a better idea of what the quotient is like in the
general case, observe that ${\Bbb C}^{\Sigma(1)} -Z(\Sigma)$ is the
union of affine open sets $U_\sigma = \{\hat x_\sigma \ne 0\}$ for
$\sigma \in \Sigma$.  Then the categorical quotient $U_\sigma/G$ is
determined by the ring of invariants of $G$ acting on the coordinate
ring of $U_\sigma$.  In \cite{Cox2}, this ring of invariants is
identified with the semigroup algebra ${\Bbb C}[\sigma^\vee\cap M]$,
which explains why the quotients $U_\sigma/G$ patch together to give
$X$.

	An interesting aspect of this construction is that it was
found by several people independently at approximately the same time.
Audin \cite{Aud} (following ideas of Delzant \cite{Del} and Kirwan
\cite{Kir}) described the construction in the context of symplectic
geometry, while Musson \cite{Mus} was studying differential operators
on toric varieties, and Cox \cite{Cox2} was more interested in the
algebraic aspects of the situation.  This result was also discovered
by Batyrev \cite{Bat3} and Fine \cite{Fin}.

	It is possible to develop the entire theory of toric
varieties using Theorem 2.1 as the {\it definition\/} of toric
variety.
\enddemo

	\subhead The exceptional set\endsubhead In $X =
\big({\Bbb C}^{\Sigma(1)} -Z(\Sigma)\big)/G$, note that the group $G$
depends only on the $1$-dimensional cones of $\Sigma$, while the
exceptional set $Z(\Sigma)$ depends on the full fan.  The
combinatorics of $Z(\Sigma)$ are quite interesting.  For projective
space, this set is very small, but it is usually bigger as the
following result of \cite{BC{\rm, Sect.~2}} shows.

\proclaim{Proposition 2.2} Let $X$ be a $n$-dimensional complete
simplicial toric variety with fan $\Sigma$.  Then either
\roster
\item $2 \le \roman{codim}\, Z(\Sigma) \le \lfloor n/2\rfloor +1$, or
\item $Z(\Sigma) = \{0\}$ and $X$ is a finite quotient of a weighted
projective space. 
\endroster
\endproclaim

	When $X$ is simplicial, the coordinate subspaces making up
$Z(\Sigma)$ can be described in terms of Batyrev's notion of a {\it
primitive collection\/}, which is a subset $\Cal{P} \subset \Sigma(1)$
with the property that $\Cal{P}$ is not the set of generators of a
cone in $\Sigma$ while every proper subset of $\Cal{P}$ is.  Then the
decomposition of $Z(\Sigma)$ into irreducible components is given by
$$
Z(\Sigma) = \bigcup_{\Cal{P}} \bold{A}(\Cal{P}),
$$
where $\bold{A}(\Cal{P})$ is the coordinate subspace determined by
$x_\rho$ for $\rho \in \Cal{P}$ and the union is over all primitive
collections $\Cal{P}$.  When $X$ is smooth and complete, Batyrev has
conjectured \cite{Bat2} that the number of irreducible components of
$Z(\Sigma)$ (= the number of primitive collections) is bounded by a
constant depending only on the the Picard number of $X$.

	\subhead The homogeneous coordinate ring\endsubhead We next
explore the algebraic consequen\-ces of $X = \big({\Bbb C}^{\Sigma(1)}
-Z(\Sigma)\big)/G$.  The basic idea is that like projective space,
homogeneous coordinates allow us to define subvarieties using global
equations.  When $X$ is a simplicial toric variety, a point of $X$ has
``homogeneous coordinates'' $(t_\rho) \in {\Bbb C}^{\Sigma(1)}
-Z(\Sigma)$, which are well-defined up to the action of $G$.  The
corresponding polynomial ring is
$$
S = {\Bbb C}[x_\rho: \rho \in \Sigma(1)]
$$
with a grading induced by the action of $G$ on ${\Bbb C}^{\Sigma(1)}$.
The character group of $G$ is $A_{n-1}(X)$ and the grading on $S$ can
be viewed as defining the ``degree'' of a monomial $x^{\bold{a}} =
\Pi_\rho x_\rho^{a_\rho}$ to be $\deg(x^{\bold{a}}) = \big[\sum_\rho
a_\rho D_\rho\big] \in A_{n-1}(X)$.  With this grading, $S = {\Bbb
C}[x_\rho]$ is called the {\it homogeneous coordinate ring\/} of $X$.

\demo{Examples} 1. In the case of ${\Bbb P}^n$ or ${\Bbb
P}^n\times{\Bbb P}^m$, the coordinate rings are the classical rings of
homogeneous or bihomogeneous polynomials.  

	2. Another good example is the blow-up of ${\Bbb C}^n$ at the
origin.  We leave it as an exercise for the reader to show that the
coordinate ring is ${\Bbb C}[t,x_1,\dots,x_n]$, where $\deg(t) = -1$,
$\deg(x_i) = 1$.  Furthermore, the exceptional set in ${\Bbb C}^{n+1}$
is $Z = {\Bbb C}\times \{(0,\dots,0)\}$ and $G = {\Bbb C}^*$ acts on
${\Bbb C}^{n+1} = {\Bbb C}\times{\Bbb C}^n$ by $\mu\cdot(t,\bold{x}) =
(\mu^{-1}t,\mu \bold{x})$.  Then, given a point $(t,\bold{x}) \in
{\Bbb C}^{n+1} - Z$, we have
$$
\align
(t,\bold{x}) &\sim_G (1,t\bold{x})\quad  \,\hbox{if $t \ne 0$}\\
(0,\bold{x}) &\sim_G (0,\mu\bold{x})\quad \hbox{if $\mu \ne 0$}.
\endalign
$$
{}From this, it should be clear that $X = ({\Bbb C}^{n+1} - Z)/G$ is the
blow-up of ${\Bbb C}^{n}$ at the origin.  Notice also that the
blow-up map $X \to {\Bbb C}^{n}$ is given by $(t,\bold{x}) \mapsto
t\bold{x}$.
\enddemo

	The coordinate ring $S = {\Bbb C}[x_\rho]$ allows us to define
subvarieties of $X$ using homogeneous ideals of $S$.  The relation
between ideals and varieties is similar to what happens in ${\Bbb
P}^n$, with the ideal $B(\Sigma) = \langle \hat x_\rho : \rho \in
\Sigma(1)\rangle$ playing the role of the irrelevant ideal (see
\cite{Cox2}).  One can also define sheaves on $X$ using graded
$S$-modules.  Here is an example of how this works.

\demo{Example} Consider the graded $S$-module $\widehat\Omega^p_S$
defined as the kernel of the map  
$$
\gamma : S\otimes \Lambda^pM \longrightarrow \bigoplus_\rho
S/\langle x_\rho\rangle \otimes \Lambda^{p-1}M,
$$
where the $\rho$th component of $\gamma$ is $\gamma_\rho(f\otimes
\omega) = f \bmod x_\rho \otimes i_\rho(\omega)$ and
$i_\rho(\omega)$ is interior product.  Then \cite{BC{\rm,
Sect.~8}} shows that the sheaf corresponding to $\widehat\Omega^p_S$
is the sheaf of Zariski $p$-forms $\widehat\Omega_X^p$ on $X$.
Furthermore, given $\alpha \in A_{n-1}(X)$, we can define a shifted
module $\widehat\Omega^p_S(\alpha)$ in the usual way.  This gives a
sheaf $\widehat\Omega^p_X(\alpha)$ with global sections
$$
H^0(X,\widehat\Omega^p_X(\alpha)) \simeq (\widehat\Omega^p_S)_\alpha
$$
(where the subscript refers to the graded piece in degree $\alpha$).  
\enddemo

	If $\Cal{L}$ is a line bundle (or, more generally, a rank one
torsion-free reflexive sheaf) on $X$, then we get $\alpha = [\Cal{L}]
\in A_{n-1}(X)$, and one can prove that
$$
H^0(X,\Cal{L}) \simeq S_\alpha\tag2.4
$$
(see \cite{Cox2}).  When $\Cal{L} = \Cal{O}_X(D)$ and $X$ is complete,
\thetag{1.2} then shows that $\dim S_\alpha$ is the number $l(\Delta)$
of integer points in the polytope $\Delta_D$.  From \thetag{2.4} we
also obtain a ring isomorphism
$$
\bigoplus_{k=0}^\infty H^0(X,\Cal{L}^{\otimes k}) \simeq
\bigoplus_{k=0}^\infty S_{k\alpha} \subset S.\tag2.5
$$
In particular, if $X$ is projective, the ``coordinate ring'' $S$
contains the coordinate rings (in the usual sense) of {\it all\/}
possible projective embeddings of $X$.

	\subhead Applications\endsubhead There have been several recent
applications of global coordinates for toric varieties.  We've already
mentioned \cite{Mus}, which studies differential operators on toric
varieties.  In \cite{Per1} and \cite{Per2}, homogeneous coordinates on
a smooth toric variety $X$ and the Euler sequence
$$
0 \longrightarrow \Omega^1_X \longrightarrow
{\textstyle{\bigoplus_\rho}} \Cal{O}_X(-D_\rho) \longrightarrow
A_{n-1}(X)\otimes\Cal{O}_X
\longrightarrow 0
$$
(see \cite{BC} and \cite{Jac}) are used to compute the principal parts
of line bundles and to study highly inflected toric varieties in
dimensions $\le 3$.

	Using homogeneous coordinates, maps to a toric variety can be
studied in much the same way one describes maps to projective
space---see, for example, \cite{Cox1}, \cite{Gue} and
\cite{Jac}.  Homogeneous coordinates were also used in \cite{Cox2} to
show that Demazure's results on the automorphism group of a smooth
complete toric variety (see \cite{Oda1{\rm, Sect.~3.4}} for a
description) remain valid in the simplicial case.

	Further applications of homogeneous coordinates will be given
in \S\S4 and 6, and in \S10, we will also see how homogeneous
coordinates are used in mirror symmetry.

\head \S3. The K\"ahler cone\endhead

	The K\"ahler classes of a smooth projective variety $X =
X_\Sigma$ form an open cone in $H^{1,1}(X,{\Bbb R})$ called the {\it
K\"ahler cone\/}.  This cone can be complicated in general, but it is
pleasantly simple when $X$ is a smooth projective toric variety.

	In this situation, the cohomology classes $[D_\rho]$ span
$H^2(X,{\Bbb R}) = H^{1,1}(X,{\Bbb R}) = A_{n-1}(X)\otimes{\Bbb R}$.
Since $X$ is smooth, a class $\bold{a} = [\sum_\rho a_\rho D_\rho]$
(with $a_\rho \in {\Bbb R}$) has a {\it support function\/} $\psi :
N_{{\Bbb R}} \to {\Bbb R}$ with the property that for each $\sigma \in
\Sigma$, there is $m_\sigma \in M_{{\Bbb R}}$ with the property that
$\psi(\rho) = \langle m_\sigma,\rho\rangle = -a_\rho$.  This is
similar to the support functions considered in \S2 except that $\psi$
is only well defined up to a linear function on $N_{{\Bbb R}}$ (this
follows from \thetag{2.2}) and the $m_\sigma$ need not be integral.
Then we say that $\bold{a}$ is {\it convex\/} if $\psi$ is a convex
function on $N_{{\Bbb R}}$.  The convex classes form a cone
$\roman{cpl}(\Sigma) \subset A_{n-1}(X)\otimes{\Bbb R}$ which has the
following nice structure.

\proclaim{Proposition 3.1} If $X$ is a simplicial projective toric
variety, then $\roman{cpl}(\Sigma) \subset A_{n-1}(X)\otimes{\Bbb R} =
H^{1,1}(X,{\Bbb R})$ is a strongly convex polyhedral cone with
nonempty interior in $H^{1,1}(X,{\Bbb R})$.  Furthermore, the interior
of this cone is precisely the K\"ahler cone of $X$.
\endproclaim

	As observed in \cite{Bat3}, the first part of the proposition
follows from \cite{OP} (see also \cite{Rei2}), and the second part
follows easily in the smooth case.  When $X$ is simplicial,
\cite{AGM2} gives a careful definition the K\"ahler cone and shows
that Proposition 3.1 continues to hold in this case.  We will see
later that this proposition has implications for both symplectic
geometry and mirror symmetry.

	Support functions $\psi$ corresponding to K\"ahler classes are
{\it strictly convex\/} and are described in \cite{Bat3} using the
primitive collections from \S2 as follows .

\proclaim{Proposition 3.2} If $X$ is a simplicial projective toric
variety, then a support function $\psi$ coming from $[\sum_\rho a_\rho
D_\rho] \in A_{n-1}(X)\otimes{\Bbb R}$ is strictly convex if and only
if for every primitive collection $\Cal{P} = \{\rho_1,\dots,\rho_k\}$,
we have
$$
\psi(\rho_1 + \cdots + \rho_k) > \psi(\rho_1) + \cdots + \psi(\rho_k).
$$
\endproclaim

	The dual of the K\"ahler cone is the {\it Mori cone\/} of
effective $1$-cycles modulo numerical equivalence.  Then Proposition
3.2 can be interpreted as describing generators for the Mori cone (see
\cite{Bat2}, \cite{OP} and \cite{Rei2} for more on the Mori cone).

\demo{Example} Consider fans $\Sigma$ in ${\Bbb R}^3$ whose
$1$-dimensional cone generators are
$$
e_0 = (0,0,-2),\ e_1 = (1,1,1),\ e_2 = (1,-1,1),\ e_3 = (-1,-1,1),\ e_4 =
(-1,1,1).
$$
Think of $e_1,e_2,e_3,e_4$ as the upper vertices of a cube and $e_0$
as lying on the negative $z$-axis.  We will use the integer lattice
generated by $e_1,e_2,e_3$.  Note that $e_0 = -e_1-e_3$ and $e_4 = e_1
- e_2 + e_3$.

	There are several ways to get a complete fan from these
generators.  For example, the cones $\sigma_{1234}$, $\sigma_{012}$,
$\sigma_{023}$, $\sigma_{034}$ and $\sigma_{041}$ (where
$\sigma_{1234}$ is the cone with generators $e_1,e_2,e_3,e_4$,
etc.)~and their faces determine a singular fan $\Sigma$.  But if we
subdivide $\sigma_{1234}$ into $\sigma_{123}$ and $\sigma_{341}$, then
we get a smooth fan $\Sigma_1$.  Similarly, we can subdivide
$\sigma_{1234}$ into $\sigma_{124}$ and $\sigma_{234}$ to get another
smooth fan $\Sigma_2$.  The toric varieties corresponding to $\Sigma$,
$\Sigma_1$ and $\Sigma_2$ will be denoted $X$, $X_1$ and $X_2$
respectively.

	The primitive collections for $\Sigma_1$ are $\{e_2,e_4\}$ and
$\{e_0,e_1,e_3\}$.  Hence, using Proposition 3.2, we see that a
support function $\psi$ is strictly convex if and only if
$$
\align
\psi(e_2+e_4) &> \psi(e_2) + \psi(e_4) \\
\psi(e_0+e_1+e_3) &> \psi(e_0) + \psi(e_1) + \psi(e_3).
\endalign
$$
If we let $\psi(e_i) = -a_i$ and use the relations $e_2+e_4 = e_1+e_3$
and $e_0+e_1+e_3 = 0$, these inequalities are equivalent to
$$
\aligned
a_2 + a_4  &> a_1 + a_3\\
a_0 + a_1 + a_3&> 0.
\endaligned\tag3.1
$$

	To determine the K\"ahler cone of $X_1$, we have to interpret
\thetag{3.1} in terms of the Chow group $A_2(X_1)\otimes{\Bbb R}$.
However, for $X_1$, the exact sequence \thetag{2.2} can be written
$$
0 \longrightarrow {\Bbb Z}^3 \overset\alpha\to\longrightarrow {\Bbb
Z}^5 \overset\beta\to\longrightarrow {\Bbb Z}^2 \longrightarrow
0,\tag3.2
$$
where $\beta$ maps $(a_0,a_1,a_2,a_3,a_4)$ to $(s,t) =
(a_0+a_1+a_3,a_0+a_2+a_4)$.  Using $s,t$ as coordinates on
$A_2(X_1)\otimes{\Bbb R} \simeq {\Bbb R}^2$, the inequalities
\thetag{3.1} can be written $t > s > 0$.  By Proposition 3.1, this is
the K\"ahler cone of $X_1$.  Thus $\roman{cpl}(\Sigma_1)$ is
$t \ge s \ge 0$.

	We now turn our attention to $X_2$.  The primitive collections
for $\Sigma_2$ are $\{e_1,e_3\}$ and $\{e_0,e_2,e_4\}$, which gives
inequalites similar to \thetag{3.1} (just interchange $a_1,a_3$ with
$a_2,a_4$).  Since \thetag{3.2} depends only on the $1$-dimensional
cones of a fan, we see that $A_2(X_1)\otimes{\Bbb R}$ is the same
${\Bbb R}^2$ with the same coordinates $s,t$.  The only difference is
that the K\"ahler cone of $X_2$ is given by $s > t > 0$ and
$\roman{cpl}(\Sigma_2)$ is  $s \ge t \ge 0$

	Hence the first quadrant in ${\Bbb R}^2$ is divided into cones
$t \ge s \ge 0$ and $s \ge t \ge 0$ whose interiors are the K\"ahler
cones of the smooth toric varieties $X_1$ and $X_2$.  Also, the
ray $s = t > 0$ corresponds to ample divisors on the singular toric
variety $X$.  
\enddemo

	As we will see in \S7, this is an example of the {\it
secondary fan\/} or {\it GKZ decomposition\/}.  We should also mention
that $X_1$ and $X_2$ are related by a {\it flop\/} (see \cite{Rei2}).

\head \S4. Symplectic geometry\endhead

	Besides the quotient ${\Bbb P}^n = ({\Bbb
C}^{n+1}-\{0\})/{\Bbb C}^*$ considered in \S2, there is the related
quotient
$$ 
{\Bbb P}^n = S^{2n+1}/S^1 
$$
where $S^1 \subset {\Bbb C}^*$ acts on the unit sphere $S^{2n+1}
\subset {\Bbb C}^{n+1}$ in the usual way.  In this section, we will
use {\it symplectic reduction} to generalize this construction to
simplicial projective toric varieties.

	\subhead The construction\endsubhead Let $X = X_\Sigma$ be the
toric variety determined by a fan $\Sigma$ in $N_{{\Bbb R}} \simeq
{\Bbb R}^n$.  We will assume that $X$ is simplicial and projective.
To simplify notation, let $r = |\Sigma(1)|$.  As in \S2, we have the
group $G = \roman{Hom}_{{\Bbb Z}}(A_{n-1}(X),{\Bbb C}^*)$.  The
maximal compact subgroup of $G$ is
$$
G_{{\Bbb R}} = \roman{Hom}_{{\Bbb Z}}(A_{n-1}(X),S^1).\tag4.1
$$
The inclusion $G \subset ({\Bbb C}^*)^r$ gives an action of $G_{{\Bbb
R}}$ on ${\Bbb C}^r$.

	Now consider the map (called the {\it moment map\/}) 
$$
\mu_\Sigma : {\Bbb C}^r \overset\mu\to\longrightarrow {\Bbb R}^r
@> {\beta_{{\Bbb R}}} >> A_{n-1}(X)\otimes{\Bbb R}, \tag4.2
$$
where $\mu$ defined by $\mu(z_1,\dots,z_r) = \frac{1}{2}
(|z_1|^2,\dots,|z_r|^2)$ and $\beta_{{\Bbb R}}$ comes from
the exact sequence
$$
0 \longrightarrow M_{{\Bbb R}} \longrightarrow {\Bbb R}^r @>
{\beta_{{\Bbb R}}} >> A_{n-1}(X)\otimes{\Bbb R} \longrightarrow 0
$$
obtained by tensoring \thetag{2.2} with ${\Bbb R}$.  Note that
$\mu_\Sigma$ is constant on $G_{{\Bbb R}}$-orbits.

	Since $X$ is projective and simplicial,
$A_{n-1}(X)\otimes{\Bbb R} \simeq H^2(X,{\Bbb R}) \simeq
H^{1,1}(X,{\Bbb R})$.  Recall from \S3 that the {\it K\"ahler cone\/}
in $H^{1,1}(X,{\Bbb R})$ consists of all possible K\"ahler classes on
$X$.  Under these isomorphisms, we get a cone in
$A_{n-1}(X)\otimes{\Bbb R}$, also called the K\"ahler cone.  Then we
modify the construction $X = \big({\Bbb C}^r - Z(\Sigma)\big)/G$ from
\S2 as follows.

\proclaim{Theorem 4.1} Let $X = X_\Sigma$ be a projective simplicial
toric variety, and assume that $\bold{a} \in A_{n-1}(X)\otimes{\Bbb
R}$ is in the K\"ahler cone.  Then $\mu_\Sigma^{-1}(\bold{a}) \subset
{\Bbb C}^r - Z(\Sigma)$, and the natural map
$$
\mu_\Sigma^{-1}(\bold{a})/G_{{\Bbb R}} \longrightarrow \big({\Bbb C}^r
- Z(\Sigma)\big)/G = X
$$
is a diffeomorphism which preserves the class of the symplectic form
(to be explained below).
\endproclaim

\demo{Proof} When $X$ is smooth, a proof that we have a diffeomorphism
can be found in Guillemin's recent book \cite{Gui{\rm, Appendix 1}},
and the statement about the class of the symplectic from follows from
equation (1.6) of \cite{Gui{\rm, Appendix 2}} (the $\lambda_i$ in
\cite{Gui} are $-a_i$ in our notation).  This proof can be modified to
work in the simplicial case.  One can also use results in \cite{Kir}
to prove the theorem.  The version in \cite{Aud} is somewhat
incomplete since the K\"ahler cone is not mentioned.\qed
\enddemo

\demo{Examples} 1. When $X = {\Bbb P}^n$, the map $\beta_{{\Bbb
R}}$ in \thetag{4.2} is the map ${\Bbb R}^{n+1} \to {\Bbb R}$ defined
by $(a_0,\dots,a_n) \mapsto \sum_{i=0}^n a_i$.  Thus
$\mu_\Sigma(z_0,\dots,z_n) = \frac{1}{2}\sum_{i=0}^n |z_i|^2$.  Since
the K\"ahler cone is ${\Bbb R}^+$, we see that $\mu_\Sigma^{-1}(a)$ is
a sphere for any $a > 0$.  Hence we recover the usual description of
${\Bbb P}^n$ as a quotient of the $(2n+1)$-sphere.

	2. Consider the toric varieties $X_1$ and $X_2$ from the
example in \S3.  These have the same group $G_{{\Bbb R}}$ acting on
${\Bbb C}^5$ and the same moment map $\mu_{\Sigma_1} =
\mu_{\Sigma_2}$, which we write as $\mu : {\Bbb C}^5 \to {\Bbb R}^2$.
Then $\mu$ is given by
$$
\mu(z_0,z_1,z_2,z_3,z_4) ={\textstyle{\frac{1}{2}}} (|z_0|^2 +
|z_1|^2 + |z_3|^2, |z_0|^2 + |z_2|^2 + |z_4|^2)
$$
(this follows from the description of $\beta$ in \thetag{3.2}).  Since
we determined the K\"ahler cones of $X_1$ and $X_2$ in \S3, it follows
that 
$$
\mu^{-1}(s,t)/G_{{\Bbb R}} \simeq \cases X_1, &\text{if $t > s >
0$}\\ X_2, &\text{if $s > t > 0$.}\endcases
$$
\enddemo 

	\subhead Symplectic manifolds and Hamiltonian
actions\endsubhead To understand the construction given in Theorem
4.1, one needs to discuss symplectic geometry.  A symplectic structure
on a real manifold $M$ is a closed, nondegenerate 2-form $\omega$.
The symplectic form $\omega$ converts functions into vector fields as
follows: if $f$ is a $C^\infty$ function on $M$, then there is a
unique vector field $X_f$ on $M$ with the property that $\omega(X,X_f)
= X(f)$ for any vector field $X$.  We call $X_f$ the {\it
Hamiltonian\/} of $f$.  Basic references on symplectic geometry are
\cite{Aud} and \cite{Kir}, and the reader might also want to consult
\cite{Ati1}, \cite{Ati2} and \cite{GS}.

	Now suppose that a compact connected Lie group $G_{{\Bbb R}}$
acts on $M$.  This action induces an infinitesimal action of the Lie
algebra ${\frak g}_{{\Bbb R}}$ where every $\lambda \in {\frak
g}_{{\Bbb R}}$ gives a vector field $X_\lambda$ on $M$.  Then the
action is {\it Hamiltonian\/} if:
\roster
\item The symplectic form $\omega$ is invariant under the group
action.
\item For each $\lambda \in {\frak g}_{{\Bbb R}}$, the vector field
$X_\lambda$ is Hamiltonian (i.e., is the Hamiltonian vector field of
some $C^\infty$ function on $M$).
\endroster
A basic property of a Hamiltonian action is that it has a {\it
moment map\/}
$$
\mu : M \longrightarrow {\frak g}_{{\Bbb R}}^*,
$$
which has the property that for every $\lambda \in {\frak
g}_{{\Bbb R}}$, the vector field $X_\lambda$ is the Hamiltonian of the
function $\lambda\circ\mu : M \to {\Bbb R}$.

\demo{Examples} 1. The most basic example is ${\Bbb C}^r$ endowed with
the symplectic form $\omega = \sum_{j=1}^r dx_j\wedge dy_j$, where
$z_j = x_j + iy_j$.  It is easy to check that the natural action of
$(S^1)^r$ on ${\Bbb C}^r$ is Hamiltonian and the moment map
$$
\mu : {\Bbb C}^r \longrightarrow ({\Bbb R}^r)^*
$$
is defined by $\mu(z_1,\dots,z_r) = \frac{1}{2}(|z_1|^2,\dots,|z_r|^2)$
(this uses the basis of $({\Bbb R}^r)^*$ dual to the standard basis of
the Lie algebra ${\Bbb R}^r$ of $(S^1)^r$).

	2. When $G_{{\Bbb R}}$ is the group coming from a toric
variety $X$ as described in \thetag{4.1}, then the action of $G_{{\Bbb
R}}$ on ${\Bbb C}^r$ is Hamiltonian and the moment map is {\it
exactly\/} as described in \thetag{4.2} provided we identify ${\Bbb
R}^r$ with it dual $({\Bbb R}^r)^*$.  Thus Theorem 4.1 tells us how
to construct a toric variety using the moment map of a Hamiltonian
action.
\enddemo

	\subhead Symplectic reduction\endsubhead Given a symplectic
manifold $M$ with a Hamiltonian action of $G_{{\Bbb R}}$ and moment
map $\mu : M \to {\frak g}_{{\Bbb R}}^*$, we can use this data to
construct other symplectic manifolds by the process of {\it symplectic
reduction\/}.  If $\bold{a} \in {\frak g}_{{\Bbb R}}$ is a regular
value of $\mu$, then $\mu^{-1}(\bold{a})$ is a manifold, but the
restriction of the symplectic form $\omega$ to $\mu^{-1}(\bold{a})$
will fail to be symplectic (it won't be nondegenerate).  However, if
$G_{{\Bbb R}}$ acts freely on $\mu^{-1}(\bold{a})$, then the
restriction of $\omega$ descends to the quotient
$\mu^{-1}(\bold{a})/G_{{\Bbb R}}$ as a symplectic form.  This is what
we mean by symplectic reduction.

	If we look back at Theorem 4.1, we see that the basic
assertion of the theorem is that we can construct smooth projective
toric varieties by symplectic reduction.  Furthermore, we can now
explain what it means for the diffeomorphism
$$
\mu_\Sigma^{-1}(\bold{a})/G_{{\Bbb R}} \longrightarrow \big({\Bbb C}^r
- Z(\Sigma)\big)/G = X
$$
to preserve the class of the symplectic form: the symplectic reduction
$\mu_\Sigma^{-1}(\bold{a})/G_{{\Bbb R}}$ has a natural symplectic
structure coming from $\omega$ on ${\Bbb C}^r$, and $X$ has a
symplectic structure coming from any K\"ahler form whose cohomology
class is $\bold{a} \in A_{n-1}(X)\otimes{\Bbb R}$.  The above
diffeomorphism need not map one symplectic form to the other, but it
does preserve their cohomology classes.

	\subhead Delzant polytopes\endsubhead In addition to the
Hamiltonian action used to construct a smooth projective toric variety
$X$, we also have an action of the real torus $T_{{\Bbb R}} = (S^1)^n$
on $X$.  It is well-known that if we give $X$ the symplectic structure
coming from an ample divisor $D = \sum_\rho a_\rho D_\rho$ (so $a_\rho
\in {\Bbb Z}$), then the action of $T_{{\Bbb R}}$ is Hamiltonian, and
the moment map
$$
\mu_X : X \longrightarrow M_{{\Bbb R}}
$$
(note that $M_{{\Bbb R}} = N^*_{{\Bbb R}}$ and $N_{{\Bbb R}}$ is the
Lie algebra of $T_{{\Bbb R}}$) can be described explicitly (see
\cite{Ful{\rm, Sect.~4.2}} or \cite{Oda1{\rm, Sect.~2.4}}).
Furthermore, the image of the moment map is, up to translation, {\it
precisely\/} the polytope $\Delta_D$ defined in \thetag{1.1}, and the
induced map
$$
X/T_{{\Bbb R}} \longrightarrow \Delta_D
$$
is a homeomorphism.

	With the publication of Delzant's thesis \cite{Del}, it became
possible to view the above results in a broader context.  Namely, given
a smooth projective toric $X$ and a class $[D] = [\sum_\rho a_\rho
D_\rho]$ in its K\"ahler cone (so $a_\rho \in {\Bbb R}$), we get the
polytope 
$$
\Delta_D = \{m \in M_{{\Bbb R}}: \langle m,\rho\rangle \ge -a_\rho\}
\subset M_{{\Bbb R}}.
$$ 
This polytope is no longer integral, but the normals to its facets are
integral, and since $X$ is smooth, the normals to the facets meeting
at each vertex of $\Delta_D$ form a basis of the lattice $N$ (this is
just another way of saying that the normal fan of $\Delta_D$ is the
fan of $X$).  Such a polytope is called a {\it Delzant polytope\/}.

\demo{Example} Consider the toric varieties $X_1$ and $X_2$ from
the example in \S3.   we leave it for the reader to show that, up to
translation, the Delzant polytope in ${\Bbb R}^3$ corresponding to $D
= a_0D_0 + \cdots + a_4D_4$ is defined by the inequalities
$$
x,y,z \ge 0,\quad x+z \le s,\quad x-y+z \ge s-t,
$$
where as usual $s = a_0+a_1+a_3$ and $t = a_0+a_2+a_4$.  For $t > s >
0$, this gives the Delzant polytope for $X_1$, and for $s > t > 0$ we
get the polytope for $X_2$.  It is a good exercise to draw these
polytopes so you can see what happens when $s=t$.  (For the best
picture, have the first octant face away from you and let the
$xz$-plane be horizontal.)
\enddemo

	If we give the toric variety $X$ a symplectic structure whose
class lies in $[D]$, then the action of $T_{{\Bbb R}}$ is still
Hamiltonian and the Delzant polytope $\Delta_D$ is again the image of
the moment map.  Conversely, given {\it any\/} Delzant polyotope
$\Delta$, one can construct a smooth projective toric variety
$X_\Delta$ with $\Delta$ as the image of the moment map (see
\cite{Gui}).

	What is more remarkable is Delzant's purely symplectic
characterization of smooth projective toric varieties (see
\cite{Del}).

\proclaim{Theorem 4.2} Let $X$ be a real $2n$-dimensional compact
connected symplectic manifold with an effective Hamiltonian action of
$(S^1)^n$.  Then the image of the moment map is a Delzant polytope
$\Delta$ and $X$ is diffeomorphic (as a Hamiltonian $(S^1)^n$-space)
to the smooth projective toric variety $X_\Delta$ determined by
$\Delta$.
\endproclaim

	One way to understand this theorem is to observe that for an
effective Hamiltonian action of $(S^1)^m$ on a connected symplectic
manifold $M$, we always have $2m \le \dim M$.  Thus Theorem 4.2
characterizes what happens when $M$ is compact and the dimension of
the torus $(S^1)^m$ is as large as possible.

\head \S5. Torus coordinates and toric ideals\endhead 

	Besides the homogeneous coordinates of \S2, toric varieties
have intrinsic coordinates living on the torus.  This is because a
basis $e_1,\dots,e_n$ of $M$ induces an isomorphism $T =
N\otimes_{{\Bbb Z}}{\Bbb C}^* \simeq ({\Bbb C}^*)^n$, giving
coordinates $t_1,\dots,t_n$ on $T$.  Then, for $m = \sum_{i=1}^n a_i
e_i \in M$, the character $\chi^m$ from \S1 is the Laurent monomial
$\bold{t}^m = \Pi_{i=1}^n t_i^{a_i}$, and the coordinate ring of $T$
is ${\Bbb C}[t_1^{\pm1},\dots,t_n^{\pm1}]$.  These coordinates don't
extend to the whole toric variety, but they are still very useful.  


	\subhead The toric variety of a polytope\endsubhead Given an
$n$-dimensional integral convex polytope $\Delta \subset M_{{\Bbb
R}}$, we get the toric variety $X_\Delta$ and an ample divisor $D$ as
mentioned in \S1.  Then (1.2) can be written as
$$
H^0(X_\Delta,\Cal{O}_X(D)) = \bigoplus_{m\in \Delta\cap M} {\Bbb
C}\cdot\bold{t}^m = L(\Delta),\tag{5.1}
$$ 
so that we can think of the global sections of $D$ in terms of Laurent
polynomials.  

	The most concrete way to get $D$ from $\Delta$ is to represent
$\Delta \subset M_{{\Bbb R}}$ by inequalities $\langle m,\rho\rangle
\ge -a_\rho$ (so the $\rho$'s are normals to the facets of $\Delta$),
and then $D = \sum_\rho a_\rho D_\rho$.  From a more sophisticated
point of view, $D$ is the divisor associated to the support function
$\psi : N_{{\Bbb R}} \to {\Bbb R}$ defined by
$$
\psi(u) = \min\{\langle m,u\rangle : m \in \Delta\}.
$$

	Since $D$ is ample, the global sections of some multiple $kD$
give a projective embedding of $X_\Delta$.  We can use this to
construct $X_\Delta$ as follows.  Let $l(k\Delta) = \dim L(k\Delta)$
be the number of integer points of $k\Delta$, and consider the map
$$
\Psi : ({\Bbb C}^*)^n \longrightarrow {\Bbb P}^{l(k\Delta) - 1}
$$
defined by 
$$
\Psi(t_1,\dots,t_n) =
(\bold{t}^{m_1},\dots,\bold{t}^{m_{l(k\Delta)}})\tag5.2
$$
where $k\Delta\cap M = \{m_1,\dots,m_{l(k\Delta)}\}$.  Then, since
$\Psi$ extends to an embedding of $X_\Delta$, it is clear that
$X_\Delta$ is the closure of $\Psi(({\Bbb C}^*)^n)$ in ${\Bbb
P}^{l(k\Delta) - 1}$.  Later in this section we will use this approach
to define non-normal toric varieties.

	A more algebraic method of constructing $X_\Delta$ is due to
Batyrev \cite{Bat4}.  Given $\Delta$, consider the cone over
$\Delta\times\{1\} \subset M_{{\Bbb R}}\oplus{\Bbb R}$ (in the
terminology of \cite{BB1}, this is a {\it Gorenstein cone\/}).  The
integer points of the cone give a semigroup algebra $S_\Delta$.  Since
$(m,k) \in M\oplus{\Bbb Z}$ is in the cone if and only if $m \in
k\Delta$, $S_\Delta$ is the subring of ${\Bbb
C}[t_0,t_1^{\pm1},\dots,t_n^{\pm1}]$ spanned by Laurent monomials
$t_0^k\bold{t}^m$ with $k \ge 0$ and $m \in k\Delta$.  This ring can
be graded by setting $\deg(t_0^k\bold{t}^m) = k$, and one can show
that
$$
X_\Delta = \roman{Proj}(S_\Delta).
$$
Since $S_\Delta$ is the coordinate ring of an affine toric variety, it
is Cohen-Macaulay and hence $X_\Delta$ is arithmetically
Cohen-Macaulay.  This ring will play an important role in the next
section.

	We can also relate $S_\Delta$ to the coordinate ring $S$ of
\S2.  If $\alpha = [D] \in A_{n-1}(X_\Delta)$ is the class of $D =
\sum_\rho a_\rho D_\rho$, then by \cite{BC}, we can define a ring
isomorphism
$$
S_\Delta \simeq \bigoplus_{k=0}^\infty S_{k\alpha} \subset S\tag5.3
$$
by sending the Laurent monomial $t_0^k\bold{t}^m$ to the monomial
$\Pi_\rho x_\rho^{ka_\rho +\langle m,\rho\rangle}$.  This is a special
case of the isomorphism \thetag{2.5}.

	\subhead Newton polytopes\endsubhead There are many situations
where instead of a polytope, the initial data is a Laurent polynomial
corresponding to a finite set $\Cal{A} \subset M = {\Bbb Z}^n$ of
exponents, which we write as
$$
f = \sum_{m \in \Cal{A}} c_m \bold{t}^m,\quad c_m \ne 0.
$$
The convex hull $\Delta = \roman{Conv}(\Cal{A})$ is called the {\it
Newton polytope\/} of $f$ and is used in many contexts (see, for
example, \cite{AVG}, \cite{Dan3}, \cite{DL}, \cite{GKZ1}, \cite{Kho},
\cite{Kus}, \cite{McD}, \cite{Var}).  This polytope might not be
$n$-dimensional, but if we use the lattice ${\Bbb Z}\Cal{A}$ generated
by $\Cal{A}$, then we get a toric variety $X_\Delta$ of dimension $=
\roman{rank}({\Bbb Z}\Cal{A})$ and, as in \thetag{5.1}, an ample line
bundle on $X_\Delta$ whose global sections are $L(\Delta) =
\oplus_{m\in\Delta\cap{\Bbb Z}\Cal{A}} {\Bbb C}\cdot \bold{t}^m$.  In
particular, our Laurent polynomial $f$ is a global section.

	A variant of this is that sometimes one is given Laurent
polynomials $f_1,\dots,f_s$ corresponding to possibly different sets
of exponents $\Cal{A}_1,\dots,\Cal{A}_s \subset {\Bbb Z}^n$.  In this
situation, we get polytopes $\Delta_i = \roman{Conv}(\Cal{A}_i)$, and
we let
$$
\Delta = \Delta_1 + \cdots + \Delta_s\tag5.4
$$ 
be their {\it Minkowski sum\/}.  Also let $L(\Delta_i) =
\oplus_{m\in\Delta_i\cap{\Bbb Z}\Cal{A}} {\Bbb C}\cdot \bold{t}^m$,
where $\Cal{A} = \Cal{A}_1\cup\cdots\cup\Cal{A}_s$.  From Proposition
2.4 of \cite{BB3}, we get the following result which is useful when
studying complete intersections in toric varieties.

\proclaim{Proposition 5.1} Given $f_1,\dots,f_s$ and $\Delta$ as
above, the toric variety $X_\Delta$ has divisors $D_i$ whose global
sections are $L(\Delta_i)$, and $\Cal{O}_{X_\Delta}(D_i)$ is generated
by these sections.  In particular, each $f_i$ is a global section of
$\Cal{O}_{X_\Delta}(D_i)$.
\endproclaim

	\subhead Non-normal toric varieties\endsubhead We can
generalize the construction \thetag{5.2} as follows.  Given a finite
set of exponents $\Cal{A} = \{m_1,\dots,m_\ell\} \subset {\Bbb Z}^n$,
we get a map
$$
\Psi : ({\Bbb C}^*)^n \longrightarrow {\Bbb C}^\ell
$$
defined by
$$
\Psi(t_1,\dots,t_n) =
(\bold{t}^{m_1},\dots,\bold{t}^{m_\ell}).\tag5.5
$$
The Zariski closure of $\Psi(({\Bbb C}^*)^n) \subset {\Bbb C}^\ell$ is
an affine variety denoted ${}^a\!X_{\Cal{A}}$ (the superscript ${}^a$
refers to affine).

	If $d = \roman{rank}({\Bbb Z}\Cal{A})$, then one can show that
$T = {}^a\!X_{\Cal{A}}\cap({\Bbb C}^*)^\ell$ is isomorphic to $({\Bbb
C}^*)^d$ and is Zariski open in ${}^a\!X_{\Cal{A}}$.  Furthermore, the
natural action of $T$ on itself extends to an action on
${}^a\!X_{\Cal{A}}$ (see \cite{Stu2{\rm, Lemma 13.4}}).  This sounds
like the definition of toric variety, but {\it normality\/} is
missing.  Nevertheless, we will refer to ${}^a\!X_{\Cal{A}}$ as a
toric variety, and we will see in this section and in \S9 that these
toric varieties are very useful.  Basic references for non-normal
toric varieties are \cite{GKZ1} and \cite{Stu2}.

\demo{Example} When $\Cal{A} = \{2,3\} \subset {\Bbb Z}$, we get the
map $\Psi(t) = (t^2,t^3)$, which leads to the cuspidal cubic $y^2 =
x^3$ in ${\Bbb C}$.  This is clearly a non-normal toric variety.
\enddemo

	The following proposition explains how ${}^a\!X_{\Cal{A}}$
relates to the usual kind of affine toric variety.  See \cite{Stu2}
for the proof.

\proclaim{Proposition 5.2} The normalization of ${}^a\!X_{\Cal{A}}$ is
the affine toric variety $X_\sigma$, where $M = {\Bbb Z}\Cal{A}$ and
$\sigma \subset N_{{\Bbb R}}$ is dual to the convex polyhedral cone
$\roman{Cone}(\Cal{A})$ generated by $\Cal{A}$.  Hence
${}^a\!X_{\Cal{A}}$ is normal if and only if ${\Bbb N}\Cal{A} = {\Bbb
Z}\Cal{A}\cap\roman{Cone}(\Cal{A})$, where ${\Bbb N}\Cal{A}$ is the
set of all non-negative integer linear combinations of $\Cal{A}$.
\endproclaim

	We next consider projective toric varieties defined using
$\Cal{A} \subset {\Bbb Z}^n$.  In practice, there are two ways of
doing this.  The first method assumes that $\Cal{A}$ lies in an affine
hyperplane in ${\Bbb Z}^n$ not passing through the origin.  In this
case, it is easy to see that ${}^a\!X_{\Cal{A}} \subset {\Bbb C}^\ell$
is defined by homogeneous polynomials.  Then ${}^a\!X_{\Cal{A}}$ is
the affine cone of a projective variety in ${\Bbb P}^{\ell-1}$ denoted
$X_{\Cal{A}}$.  One can prove that $X_{\Cal{A}}$ is a toric variety
(possibly non-normal) of dimension equal to the dimension of the
affine span of $\Cal{A}$ (= the dimension of the convex hull
$\roman{Conv}(\Cal{A})$).

	A second method for creating projective toric varieties starts
with an {\it arbitrary\/} subset $\Cal{A} = \{m_1,\dots,m_\ell\}
\subset {\Bbb Z}^n$ and considers the map
$$
\Psi(t_1,\dots,t_n) = (\bold{t}^{m_1},\dots,\bold{t}^{m_\ell}) \in
{\Bbb P}^{\ell-1} .\tag5.6
$$
Then, as in \cite{GKZ1}, we define $X_{\Cal{A}}$ to be the closure in
${\Bbb P}^{\ell-1}$ of the image of $\Psi$.  The dimension of
$X_{\Cal{A}}$ again equals the dimension of $\roman{Conv}(\Cal{A})$.

	These two approaches are related as follows.  Given $\Cal{A}$
as in the second method, $\Cal{A}\times\{1\} \subset {\Bbb Z}^{n+1}$
lies in an affine hyperplane not passing through the origin.  Then the
affine cone of $X_{\Cal{A}}$ is easily seen to be
${}^a\!X_{\Cal{A}\times\{1\}}$.  Thus $X_{\Cal{A}}$, as defined by the
second method, equals $X_{\Cal{A}\times\{1\}}$, as  defined by the first.

	The normalization of the projective toric variety
$X_{\Cal{A}}$ can be computed using Proposition 5.2.  This is a bit
delicate because of the distinction between normality and projective
normality---it is possible for $X_{\Cal{A}}$ to be normal without
${}^a\!X_{\Cal{A}}$ being so (see \cite{Har{\rm, Ex.~3.18 on
p.~23}} for an example).

\proclaim{Proposition 5.3} The normalization of $X_{\Cal{A}}$ is
$X_\Delta$, where $\Delta = \roman{Conv}(\Cal{A})$. 
\endproclaim

\demo{Proof} We regard $X_{\Cal{A}}$ as arising from \thetag{5.6}.  As
noted above, its affine cone is ${}^a\!X_{\Cal{A}\times\{1\}}$, so by
Proposition 5.2, the normalization of ${}^a\!X_{\Cal{A}\times\{1\}}$
comes from the cone over $\Cal{A}\times\{1\}$.  This equals the cone
over $\Delta\times\{1\}$, which is {\it exactly\/} the cone used in
Batyrev's construction of $X_\Delta$ earlier in this section.  Thus,
the normalization of ${}^a\!X_{\Cal{A}}$ is the affine cone of
$X_\Delta$, which implies that the map $X_\Delta \to X_{\Cal{A}}$ is
finite and birational.  The proposition now follows from Zariski's
Main Theorem.\qed
\enddemo

\demo{Example} To see how non-normal projective toric varieties can
occur, consider a complete toric variety $X$ (in 
the usual sense) with an ample divisor $D$.  If we let $\Cal{A} =
\Delta_D\cap M$ and use the basis of $H^0(X,\Cal{O}_X(D))$ given by
Laurent monomials, we get a map $\Psi : X \to {\Bbb P}^{\ell - 1}$ as
above.  Since $D$ need not be very ample, $\Psi$ need not be an
embedding.  But \thetag{5.6} shows that the image $\Psi(X)$ is
precisely the toric variety $X_{\Cal{A}}$.  It follows from
Proposition 5.3 that {\it $X$ is the normalization of the image of the
map to projective space given by an ample line bundle on $X$.}
\enddemo

	There is a nice criterion for normality which involves the
Hilbert polynomial of $X_{\Cal{A}}$ and the Ehrhart polynomial of the
polytope $\Delta = \roman{Conv}(\Cal{A})$.  By \cite{Stu2{\rm,
Ch.~13}}, the Hilbert polynomial of $X_{\Cal{A}}$ is given by
$$
\align
H_{\Cal{A}}(k) &= |\{m_{i_1} + \cdots + m_{i_k} :
m_{i_1},\dots,m_{i_k} \in \Cal{A}\}|,\quad k \gg 0,\\
\intertext{and, as usual, the Ehrhart polynomial of $\Delta$ is}
E_\Delta(k) &= |{\Bbb Z}\Cal{A}\cap k\Delta|,\quad k \ge 0
\endalign
$$
These polynomials have the same leading term, which implies that the
degree of $X_{\Cal{A}} \subset {\Bbb P}^{\ell-1}$ is the normalized
volume of $\Delta$ (see \cite{Stu2{\rm, Thm.~4.16}}).  Then we have
the following result of Sturmfels \cite{Stu2{\rm, Thm.~13.11}}.

\proclaim{Theorem 5.4} The toric variety $X_{\Cal{A}} \subset {\Bbb
P}^{\ell-1}$ is normal if and only if the Hilbert polynomial
$H_{\Cal{A}}$ equals the Ehrhart polynomial $E_{\Cal{A}}$.
\endproclaim 

	\subhead Toric ideals\endsubhead A familiar example from
algebraic geometry is the twisted cubic $(x,y,z) = (t,t^2,t^3)$ in
${\Bbb C}^3$.  This is a special case of the construction
\thetag{5.5}.  The ideal of the twisted cubic is $\langle
y-x^2,z-x^3\rangle \subset {\Bbb C}[x,y,z]$ and is our first
example of a toric ideal.  As we will soon see, the simple form of its
generators is no accident.

	Given $\Cal{A} = \{m_1,\dots,m_\ell\} \subset {\Bbb Z}^d$, the
affine toric variety ${}^a\!X_{\Cal{A}}$ is defined by an ideal
$I_{\Cal{A}} \subset {\Bbb C}[x_1,\dots,x_\ell]$, which we call a {\it
toric ideal\/}.  In terms of elimination theory, the ideal
$I_{\Cal{A}}$ arises from the equations $x_i - \bold{t}^{m_i} = 0,\
1-yt_1\cdots t_n = 0$ by eliminating $y,t_1,\dots,t_n$ (the equation
$1-yt_1\cdots t_n = 0$ guarantees that $(t_1,\dots,t_n) \in ({\Bbb
C}^*)^n$).

	A toric ideal $I_{\Cal{A}}$ is homogeneous when $\Cal{A}$
lies in an affine hyperplane missing the origin, in which case
$I_{\Cal{A}}$ defines the projective toric variety $X_{\Cal{A}}$.
Also, if $X_{\Cal{A}}$ is projectively normal, then ${\Bbb
C}[x_1,\dots,x_\ell]/I_{\Cal{A}}$ is isomorphic to the ring $S_\Delta$
of \thetag{5.3}.  So facts about toric ideals give useful information
about the coordinate rings of projective toric varieties.

	As with the twisted cubic, toric ideals are generated by {\it
binomials\/}, which are differences of monomials.  To state the
result, note that a vector $\bold{a} \in {\Bbb Z}^\ell$ can be written
$\bold{a}^+ - \bold{a}^-$, where $\bold{a}^+$ and $\bold{a}^-$ have
non-negative entries and disjoint support.

\proclaim{Lemma 5.5} If $\Cal{A} = \{m_1,\dots,m_\ell\} \subset {\Bbb
Z}^n$, then the toric ideal $I_{\Cal{A}}$ can be written
$$
I_{\Cal{A}} = \langle \bold{x}^{\bold{a}^+} - \bold{x}^{\bold{a}^-} :
\bold{a} = (a_1,\dots,a_\ell) \in {\Bbb Z}^\ell,\
\textstyle{\sum_{i=1}^\ell} a_im_i = 0\rangle.
$$
\endproclaim

	In \cite{Stu1} and \cite{Stu2}, one can find a wealth of
results about toric ideals, including facts about Gr\"obner bases and 
relations to secondary fans.  We will discuss these ideas briefly in
\S7.  Applications to enumeration, sampling, integer programming, and
primitive partition identities are given in \cite{Stu2{\rm, Ch.~5 and
6}}.

	We close this section with an unsolved problem about toric
ideals.

\proclaim{Conjecture 5.6} If $X_{\Cal{A}}$ is a smooth projectively
normal toric variety, then the toric ideal $I_{\Cal{A}}$ is generated
by quadratic binomials.
\endproclaim

	In this situation, it is known that $I_{\Cal{A}}$ is generated
by binomials of degree at most $d = \roman{rank}({\Bbb Z}\Cal{A})$
(see \cite{Stu2{\rm, Thm.~13.14}}).
	
\head \S6. Cohomology of toric hypersurfaces\endhead 

	A complete toric variety $X$, being rational, is a very
special kind of variety.  But as an ambient space, $X$ can be home to
some interesting subvarieties (we will see some Calabi-Yau examples in
\S8).  In particular, using the homogeneous coordinates of \S2 or the
torus coordinates of \S5, it is easy to describe hypersurfaces and
complete interesections in $X$.  In this section, we will study
hypersurfaces $Y \subset X$ and the associated affine hypersurfaces
$Y\cap T \subset T$, where $T$ is the torus of $X$.  We will end the
section with some remarks about complete intersections.

	\subhead Affine and projective hypersurfaces\endsubhead
Given a Laurent polynomial $f$ contained in  ${\Bbb
C}[t_1^{\pm1},\dots,t_n^{\pm1}]$, the equation $f=0$ defines an
affine hypersurface
$$
Z_f \subset T = ({\Bbb C}^*)^n.
$$
We can assume that $f \in L(\Delta)$ for some convex integral polytope
$\Delta \subset {\Bbb R}^n$ (for example, let $\Delta$ be the Newton
polytope of $f$).  We will also assume that $\Delta$ has dimension
$n$ (it is easy to reduce to this case).  Then we have a toric variety
$X_\Delta$ containing $T$, and $f \in L(\Delta)$ defines a projective
hypersurface
$$
Y_f \subset X_\Delta
$$
since $f$ can be regarded as a global section of an ample line bundle
on $X_\Delta$.  Note that $Z_f = Y_f \cap T$, though $Y_f$ need not be
the Zariski closure of $Z_f$.  The latter happens, for instance, if
$f$ is a single monomial $\bold{t}^m$.

	When $\Delta$ is the Newton polytope of $f$, there is a nice
relation between the topology of $Z_f$ and the vertices of $\Delta$.
Using the map
$$
\log : ({\Bbb C}^*)^n \longrightarrow {\Bbb R}^n,\quad (t_1,\dots,t_n)
\mapsto (\log|t_1|,\dots,\log|t_n|)\tag6.1
$$
one can show that unbounded connected components of ${\Bbb R}^n -
\log(Z_f)$ correspond to vertices of $\Delta$ and in each such
component contains a translate of the corresponding cone in the normal
fan of $\Delta$.  There is also a version of this (using the moment
map) for $Y_f$---see \cite{GKZ{\rm, Sect.~1.B and 1.C of Ch.~6}}.

	Since $Z_f$ or $Y_f$ could be very singular for an arbitrary
$f \in L(\Delta)$, we need some sort of genericity condition on $f$.
There are two conditions which are used in practice.  First, $f$ is
{\it nondegenerate\/} or {\it $\Delta$-regular\/} if $Y_f \subset
X_\Delta$ meets every torus orbit transversely, i.e., if for every
orbit $O$, $Y_f\cap O$ is smooth of codimension $1$ in $O$.  In
particular, a nondegenerate $Y_f$ misses all fixed points of $T$, and
since fixed points correspond to vertices of $\Delta$, it follows that
$\Delta$ is the Newton polytope of $f$.  By Bertini's Theorem, a
generic $f \in L(\Delta)$ is nondegenerate (see
\cite{Kho}).

	The second condition applies when $X_\Delta$ is simplicial,
which for $\Delta$ means that every vertex is contained in precisely
$n$ faces (such a polytope is called {\it simple\/}).  In this case,
$X_\Delta$ is a $V$-manifold, and $f \in L(\Delta)$ is {\it
quasi-smooth\/} if $Y_f \subset X_\Delta$ is a $V$-submanifold, as
defined in \cite{BC{\rm, Sect.~3}}.  When $X_\Delta$ is smooth, this
means that $Y_f$ is smooth.  One can also show that nondegenerate
implies quasi-smooth, though the converse is not true (the smooth
cubic $y^2z = x^3-z^3$ in ${\Bbb P}^2$ is an example).  The notion of
quasi-smooth was introduced by Danilov (see \cite{Dan2{\rm,
Sect.~14}}). 

	In the simplicial case, we can use the homogeneous coordinate
ring $S = {\Bbb C}[x_\rho]$ of \S2 to describe hypersurfaces in
$X_\Delta$.  If $\beta = [D] \in A_{n-1}(X_\Delta)$ is the class of
the line bundle $D$ determined by $\Delta$, then there is a natural
isomorphism
$$
L(\Delta) \simeq S_\beta\tag6.2
$$
described in \thetag{5.3}. Thus, $f \in L(\Delta)$ corresponds to a
homogeneous polynomial $F \in S_\beta$.  The equation $F = 0$ gives a
hypersurface in ${\Bbb C}^{\Sigma(1)}-Z(\Sigma)$ (where $\Sigma$ is
the fan of $X_\Delta$).  By Theorem 2.1, this descends to a
hypersurface $Y_F \subset X_\Delta$, and one can check that $Y_F =
Y_f$.  Furthermore, $Y_F$ is quasi-smooth if and only if the partial
derivatives $\partial F/\partial x_\rho$ have no common zeros in
${\Bbb C}^{\Sigma(1)}-Z(\Sigma)$ (see \cite{BC{\rm, Sect.~3}}).

	When the defining polynomial $f$ or $F$ is clear from context,
we will write $Z \subset T$ in the affine case and $Y \subset
X_\Delta$ in the projective case.

	\subhead Cohomology of affine hypersurfaces\endsubhead For a
nondegenerate affine hypersurface $Z = Z_f \subset T = ({\Bbb
C}^*)^n$, the cohomology groups $H^i(Z)$ (always with complex
coefficients) carry natural mixed Hodge structures.  In \cite{DK}, the
dimension $h^{p,q}(H^i(Z))$ of the $(p,q)$ Hodge component of
$\roman{Gr}^{\Cal{W}}_{p+q} H^i(Z)$ (where $\Cal{W}$ is the weight
filtration) is computed using the combinatorics of the Newton polytope
$\Delta$ of $f$.

	It suffices to compute $h^{p,q}(H^i_c(Z))$ since $H^i_c(Z)$ is
dual (as a mixed Hodge structure) to $H^{2n-2-i}(Z)$.  In fact, it is
sufficient to compute 
$$
e^{p,q}(Z) = \sum_i (-1)^i h^{p,q}(H_c^i(Z))
$$
since the Gysin map
$$
H_c^i(Z) \longrightarrow H_c^{i+2}(T)
$$
is an isomorphism of mixed Hodge structures (suitably shifted) for $i
> n-1$ and $H^i_c(Z) = 0$ for $i < n-1$ (recall that $Z$ is smooth
and affine).    

	In the special case when $\Delta$ is simplicial, computing
$e^{p,q}(Z)$ is fairly easy.  First, for a face $\Gamma \subset
\Delta$, define $\phi_i(\Gamma)$ for $0 < i \le \dim \Gamma$ by the
formulas
$$
\aligned
\phi_i(\Gamma) &= \sum_{k=1}^i (-1)^{i-k} \binom{\dim \Gamma +1}{i-k}
l^*(k\Gamma)\\
	&= (-1)^{\dim \Gamma + 1} \sum_{k=0}^{\dim \Gamma + 1 - i}
(-1)^{i+k} \binom{\dim \Gamma +1}{i+k} l(k\Gamma)
\endaligned\tag6.3
$$
where $l(k\Gamma)$ (resp.\ $l^*(k\Gamma)$) is the number of integer
points in $k\Gamma$ (resp.\ in the relative interior of $k\Gamma$).
The two representations of $\phi_i(\Gamma)$ are related to the
remarkable properties of the Ehrhart polynomial from \S5 (see
\cite{Bat4{\rm, Sect.~2}} for more details).

	When $p > q$, \cite{DK{\rm, 5.7}} gives the formula
$$
e^{p,q}(Z) = (-1)^{n+p+q} \sum_{\dim\Gamma > p} (-1)^{\dim \Gamma}
\binom{n-\dim \Gamma}{n-1-p-q}\phi_{\dim \Gamma - p}(\Gamma),
$$
where the sum is over all faces $\Gamma$ of $\Delta$ of dimension $>
p$.  Since $e^{p,q}(Z) = e^{q,p}(Z)$, this gives us everything except
$e^{p,p}(Z)$.  However, we also have the identity
$$
(-1)^{n-1}\sum_q e^{p,q}(Z) = (-1)^p\binom{n}{p+1} +
\phi_{n-p}(\Delta),\tag6.4
$$
from \cite{DK{\rm, 4.4}}, which now enables us to compute
$e^{p,p}(Z)$.  

	The paper \cite{DK} also describes an algorithm for computing
the full mixed Hodge structure of $Z_f \subset X_\Delta$ for {\it
any\/} polytope $\Delta$.  There are tables giving explicit formulas
when $1 \le \dim(\Delta) \le 4$ (see \cite{DK{\rm, 5.11}}).

	One way to represent the numbers $e^{p,q}(Z)$ is via the
$E$-polynomial
$$
E(Z;u,v) = \sum_{p,q} e^{p,q}(Z) u^p v^q.
$$
Then \thetag{6.4} describes $E(Z;u,1)$, which tells us about the Hodge
filtration.  The weight filtration, on the other hand, concerns
$E(Z;u,u)$.  This polynomial is studied in \cite{BB2} and \cite{DL}.
We should also mention that explicit formulas for $E(Z;u,v)$ (for
$\Delta$ arbitary, not just simplicial) can be found in \cite{BB2}.

	We next describe some work of Batyrev \cite{Bat4} on
representing cohomology classes of $Z \subset T$ as residues of forms.
Here, we will focus on the mixed Hodge structure of $H^i(Z)$.  The
dual of the above Gysin map is the natural map $H^k(T) \to H^k(Z)$,
which is an isomorphism for $k < n-1$ and injective for $i = n-1$.
Because of this, we define the {\it primitive cohomology\/} of $Z$ by
the exact sequence
$$
0 \longrightarrow H^{n-1}(T) \longrightarrow H^{n-1}(Z)
\longrightarrow H_0^{n-1}(Z) \longrightarrow 0.\tag6.5
$$
It follows that $H_0^{n-1}(Z)$ has a mixed Hodge structure.  The
graded pieces of the Hodge filtration will be denoted
$\roman{Gr}_{\Cal{F}}^pH_0^{n-1}(Z)$.  In \cite{Bat4{\rm, Cor.~3.14}},
Batyrev observes that \thetag{6.4} can be reformulated as
$$
\dim \roman{Gr}_{\Cal{F}}^pH_0^{n-1}(Z) = \phi_{p+1}(\Delta).\tag6.6
$$

\demo{Example} Let $\Delta$ be a convex integer polygon in
${\Bbb R}^2$.  Then $Z \subset ({\Bbb C}^*)^2$ is a curve which can be
obtained by removing $m$ points from a smooth complete curve of genus
$g$.  As is well-known, this determines the Hodge numbers of $H^1(Z)$,
and using \thetag{6.5}, we obtain $h^{1,0}(H_0^1(Z)) =
h^{0,1}(H_0^1(Z)) = g$ and $h^{1,1}(H_0^1(Z)) = m - 3$.  However,
\thetag{6.6} and \thetag{6.3} imply that
$$
\align
h^{0,1}(H_0^1(Z)) &= \phi_1(\Delta) = l^*(\Delta)\\
h^{1,0}(H_0^1(Z)) + h^{1,1}(H_0^1(Z)) &= 
	\phi_2(\Delta) = l(\Delta) - 3.
\endalign
$$
It follows that $g = l^*(\Delta)$ and $m = l(\Delta) - l^*(\Delta) =$
the number of integer points on the boundary of $\Delta$.  Thus we can
see the mixed Hodge structure of $Z$ geometrically.  This example can
also be done using the formulas of \cite{DK}. 
\enddemo

	The next step is to represent cohomology classes in
$H^{n-1}_0(Z)$ algebraically.  For this purpose, let $g = t_0
f(t_1,\dots,t_n) - 1$, which is in the ring $S_\Delta$ introduced in
\S5.  We also set
$$
g_i = t_i \frac{\partial g}{\partial t_i},\quad 0 \le i \le n.
$$
Then $g_0 = t_0 f$ and $g_i \in S_\Delta$ has degree $1$ in $S_\Delta$
because $t_0$ appears to the first power.  These polynomials generate
the graded ideal $J_{f,\Delta} = \langle g_0,\dots,g_n\rangle \subset
S_\Delta$.
The following result is proved in \cite{Bat4{\rm, Sect.~6}}.

\proclaim{Theorem 6.1} There is a natural isomorphism
$$
\roman{Gr}_{\Cal{F}}^p H_0^{n-1}(Z) \simeq
(S_\Delta/J_{f,\Delta})_{n-p}.
$$
\endproclaim

	The idea behind this isomorphism is that a polynomial
$t_0^{n-p}g(t_1,\dots,t_n) \in (S_\Delta)_{n-p}$ gives a $n$-form
$$
\frac{g}{f^{n-p}} \frac{dt_1}{t_1}\wedge\cdots\wedge \frac{dt_n}{t_n}
$$
on $T-Z$ whose residue lies in $\roman{Gr}_{\Cal{F}}^p H_0^{n-1}(Z)$.

	\subhead Cohomology of projective hypersurfaces\endsubhead
We now turn our attention to the projective hypersurface $Y = Y_f
\subset X = X_\Delta$, where as usual $f \in L(\Delta)$ is
nondegenerate.  In addition, we will assume that $X$ is simplicial.
Then the natural map
$$
H^i(X) \longrightarrow H^i(Y)
$$
is an isomorphism for $i < n-1$ and is injective for $i = n-1$.  In
analogy with the affine case, we define the {\it primitive
cohomology\/} by the exact sequence
$$
0 \longrightarrow H^{n-1}(X) \longrightarrow H^{n-1}(Y)
\longrightarrow H_0^{n-1}(Y) \longrightarrow 0.
$$
Since $X$ is simplicial, $Y$ is quasi-smooth and hence $H_0^{n-1}(Y)$
has a pure Hodge structure.  Letting $h_0^{p,n-1-p}$ denote the
dimension of the appropriate Hodge component, \cite{DK{\rm, 5.5}} can
be restated as
$$
h_0^{p,n-1-p} = (-1)^{n} \sum_{\dim\Gamma > p} (-1)^{\dim \Gamma}
\phi_{\dim \Gamma - p}(\Gamma),
$$
where $\phi_i(\Gamma)$ is as in \thetag{6.3}.  This formula implies
that $h_0^{n-1,0} = \phi_1(\Delta) = l^*(\Delta)$.

	We have already seen that $Y \subset X$ can be defined
in two ways, either using $f \in L(\Delta)$ or $F \in S_\beta$, where
$f$ and $F$ are related via \thetag{6.2}.  Then we have the following
results from \cite{Bat4} and \cite{BC}.

\proclaim{Theorem 6.2} Let $X = X_\Delta$ be simplicial and $Y = Y_f$ be
nondegenerate.  If $I^{(1)} \subset S_\Delta$ be the ideal generated by those
$t_0^k\bold{t}^m \in S_\Delta$ for which $m$ is an interior point of
$k\Delta$, and $H_f$ is its image in $S_\Delta/J_{f,\Delta}$, then
there are natural isomorphisms
$$
H^{p,n-1-p}_0(Y) \simeq \roman{Gr}_{\Cal{F}}^p
\Cal{W}_{n-1}H^{n-1}_0(Z) \simeq (H_f)_{n-p},
$$
where $Z = Z_f$ is the corresponding affine hypersurface.
\endproclaim

\proclaim{Theorem 6.3} Let $X = X_\Delta$ be simplicial and $Y = Y_F$ be
quasi-smooth.  If $J(F) = \langle \frac{\partial F}{\partial
x_\rho}\rangle \subset S$ be the Jacobian ideal of $F$, and $\beta_0
=\big [\sum_\rho D_\rho\big]$ be the anticanonical class, then there
is a natural isomorphism
$$
\,H^{p,n-1-p}_0(Y) \simeq (S/J(F))_{(n-p)\beta-\beta_0}
\phantom{\qquad\qquad\qquad}
$$
when $p \ne n/2-1$, and for $p = n/2-1$, we have an exact sequence
$$
\qquad 0 \to H^{n-2}(X)
\overset{\cup[Y]}\to\longrightarrow H^{n}(X) \to
(S/J(F))_{(n/2+1)\beta-\beta_0}\to
H^{n/2-1,n/2}_0(Y)\to 0.
$$
\endproclaim

	Theorem 6.3 generalizes a classic result of Griffiths on the
cohomology of projective hypersurfaces (see \cite{Gri} and \cite{PS}).
The basic idea is similar to Theorem~6.1: a homogeneous polynomial $G
\in S_{(n-p)\beta-\beta_0}$ gives a $n$-form on $X-Y$ whose residue
lies in $H_0^{n-1-p}(Y)$ (see \cite{BC} for a precise description).
The complication in the case $p=n/2-1$ comes from the exact sequence
$$
0 \to H^{n-2}(X)
\overset{\cup[Y]}\to\longrightarrow H^{n}(X) \to
H^n(X-Y) \to
H^{n-1}_0(Y)\to 0.
$$
When $X$ is projective space (or even a weighted projective space),
the sequence implies $H^n(X-Y) \simeq H^{n-1}_0(Y)$, but this can fail
in the toric case.  Fortunately, this only affects
$H^{n/2-1,n/2}_0(Y)$.  

\demo{Remarks} 1. While Theorem 6.1 uses the vanishing of
$H^i(X,\Cal{O}_X(Y))$ for $Y$ ample and $i > 0$, Theorem 6.3 uses the
Bott-Steenbrink-Danilov vanishing theorem:
$$
H^i(X,\widehat\Omega^p_X(Y)) = 0,\quad Y\ \roman{ample},\ i> 0.\tag6.7
$$ 
This result is stated in \cite{Dan2} and \cite{Oda2} without proof.  In
the simplicial case, a proof appeared in \cite{BC{\rm, Sect.~7}}, and
a general proof (using characteristic $p > 0$) can be found in
\cite{BFLM}.  A generalization of \thetag{6.7} is the vanishing
theorem:
$$
H^i(X,\Cal{W}_k\widehat\Omega^p_X(\log (-K))(Y)) = 0,\quad Y\
\roman{ample},\ i> 0.\tag6.8
$$ 
Here, $\Cal{W}_k\widehat\Omega^p_X(\log (-K)) = \widehat\Omega^{p-k}_X
\wedge \widehat\Omega^k_X(\log (-K))$ is the usual weight filtration
and $-K = \sum_\rho D_\rho$ is the anticanonical divisor.
When $X$ is simplicial, \thetag{6.8} was proved in \cite{BC}, but the
general case is still open.

	2. Using the isomorphism \thetag{5.3}, one can relate the
ideal $H_f \subset S_\Delta/J_{f,\Delta}$ to the ideal generated by
$\Pi_\rho x_\rho$ in $S/\langle x_\rho \partial F/\partial
x_\rho\rangle$.  This leads to a natural isomorphism
$$
H_0^{p,n-1-p}(Y) \simeq (S/J_1(F))_{(n-p)\beta-\beta_0},
$$
where $J_1(F)$ is the ideal quotient
$$
J_1(F) = \langle x_\rho \partial F/\partial x_\rho\rangle :
{\textstyle{\Pi_\rho}} x_\rho
$$
(see \cite{BC{\rm, Sect.~11}}).  For a weighted projective space,
$J(F)$ equals $J_1(F)$, but in general the relation between these
ideals is not well understood.  The ideal $J_1(F)$ arises naturally in
certain mirror symmetry contexts (see \cite{MP{\rm, 5.36}}).
\enddemo

	It can also happen that one is interested in a hypersurface $Y
\subset X = X_\Delta$ which is big and nef but not ample.  In the
toric context, big and nef mean that $Y$ corresponds to a
$n$-dimensional integer polytope $\Delta'$ which is a Minkowski
summand of $\Delta$ (i.e., $\Delta' + \Delta'' = \mu\Delta$ for some
integer polyope $\Delta''$ and $\mu \in {\Bbb Z}$---see \cite{BB3{\rm,
Sect.~2}}).  In this case, the map $H^i(X) \to H^i(Y)$ need not be an
isomorphism for $i < n-1$, and Theorems~6.2 and 6.3 can also fail.  An
example of the latter is given by the proper transform of a degree 8
hypersurface in a resolution of ${\Bbb P}(1,1,2,2,2)$.  We will say
more about this when we study Calabi-Yau hypersurfaces in \S8.

	Besides the cohomology of toric hypersurfaces, one can also
study their moduli (see \cite{AGM2}, \cite{Bat4} and \cite{BC}).  The
resulting variations of Hodge structure are closely connected with
hypergeometric functions, which will be discussed in \S9.

	\subhead Complete intersections\endsubhead Complete
intersections in toric varieties can be studied from several points of
view.  In the affine case, suppose we have $f_i \in L(\Delta_i)$ as in
Proposition 5.1.  Then we get the affine complete intersection
$$
Z_{f_1}\cap\cdots\cap Z_{f_s} \subset T.
$$
To compute the cohomology of this variety, we use the {\it Cayley
trick\/}.  Consider the toric variety $T\times{\Bbb C}^s$ with variables
$t_1,\dots,t_n,\lambda_1,\dots,\lambda_s$, and let
$$
\Cal{F} = \lambda_1 f_1 + \cdots \lambda_s f_s - 1.
$$
This gives the affine hypersurface $Z_{\Cal{F}} \subset T\times{\Bbb
C}^s$, and one obtains
$$
E(Z_{f_1}\cap\cdots\cap Z_{f_s};u,v) = (uv-1)^s - (uv)^{1-s}
E(Z_{\Cal{F}};u,v)
$$
under suitable nondegeneracy hypotheses.  This follows by considering
the projection $Z_{\Cal{F}} \to T$ (see \cite{DK{\rm, 6.2}}).  Hence
we can reduce to the hypersurface case (note, however, that
$Z_{\Cal{F}}$ is a hypersurface in $T\times{\Bbb C}^s$ instead of a
torus).

	Turning to the projective case, let $\Delta = \Delta_1 +
\cdots + \Delta_s$ be as in \thetag{5.4}, and consider the complete
intersection 
$$
Y_{f_1} \cap\cdots\cap Y_{f_s} \subset X = X_\Delta,
$$
which we assume to be nondegenerate.  Now consider the projective
bundle 
$$
{\Bbb P}(\Cal{E}) = {\Bbb
P}(\Cal{O}_X(D_1)\otimes\cdots\otimes\Cal{O}_X(D_s)), 
$$
where $D_i$ is the divisor corresponding to $\Delta_i$ in Proposition
5.1.  If $\pi : {\Bbb P}(\Cal{E}) \to X$ is the natural projection,
then $\pi_*(\Cal{O}_{{\Bbb P}(\Cal{E})}(1)) =
\Cal{O}_X(D_1)\otimes\cdots\otimes\Cal{O}_X(D_s)$.  Hence there is a
unique section of $\Cal{O}_{{\Bbb P}(\Cal{E})}(1)$ corresponding to
$(f_1,\dots,f_s)$.  This section defines a hypersurface $\Cal{Y}
\subset {\Bbb P}(\Cal{E})$, and one can show that the natural map
$$
{\Bbb P}(\Cal{E}) - \Cal{Y} \longrightarrow X - Y_{f_1} \cap\cdots\cap
Y_{f_s} 
$$
is a ${\Bbb C}^{s-1}$ bundle in the Zariski topology.  Hence
$$
H_c^i(X -  Y_{f_1} \cap\cdots\cap Y_{f_s}) \simeq
H_c^{i+2(s-1)}({\Bbb P}(\Cal{E}) - \Cal{Y}), 
$$
which is an isomorphism of mixed Hodge structures of degree
$(s-1,s-1)$.  This is discussed in more detail in \cite{BB1}.  

	From the point of view of the homogeneous coordinate ring $S$
of $X$, each $f_i$ corresponds to a homogeneous polynomial $F_i \in
S$.  If we assume that each $\Cal{O}_X(D_i)$ is ample, then it follows
from \cite{CCD} that the homogeneous coordinate ring of ${\Bbb
P}(\Cal{E})$ is $S\otimes{\Bbb C}[y_1,\dots,y_s]$, where the variables
$y_i$ have the property that
$$
F = y_1 F_1 + \cdots + y_s F_s\tag6.9
$$
is homogeneous and defines the hypersurface $\Cal{Y} \subset {\Bbb
P}(\Cal{E})$.  This explains why the above construction can be
regarded as the projective version of the Cayley trick.

	It should also possible to represent primitive cohomology
classes of $Y_{f_1} \cap\cdots\cap Y_{f_s}$ using the Jacobian ideal
of $F \in S\otimes{\Bbb C}[y_1,\dots,y_s]$.  This has not been done in
general, but the case $X = {\Bbb P}^n$ has been studied in \cite{ENV},
\cite{Konn}, \cite{LT} and \cite{Ter}.  Furthermore, \cite{Dim} and
\cite{Nag} treat this case from the toric point of view, which we now
explain.  Let $F_1,\dots,F_s$ be homogeneous polynomials (in the usual
sense) in ${\Bbb C}[x_0,\dots,x_n]$ of degrees $d_1,\dots,d_s$ at
least 2.  This gives the complete intersection $Y_{F_1}
\cap\cdots\cap Y_{F_s} \subset {\Bbb P}^n$ of dimension $n-s$.  The
homogeneous coordinate ring of ${\Bbb P}(\Cal{E})$ is $R = {\Bbb
C}[x_0,\dots,x_n,y_1,\dots,y_s]$ and is graded by ${\Bbb Z}^2$, where
$$
\deg(x_i) = (1,0),\quad \deg(y_j) = (-d_j,1).
$$
The polynomial $F$ of \thetag{6.9} has degree $\beta = (0,1)$ in $R$
and, defining primitive cohomology in the usual way, one obtains a
natural isomorphism
$$
H_0^{p,n-s-p}(Y_{F_1} \cap\cdots\cap Y_{F_s}) \simeq
(R/J(F))_{(n-p)\beta - \beta_0},\tag6.10
$$
where $\beta = \deg(F) = (0,1)$, $\beta_0 = \deg(x_0\cdots
x_ny_1\cdots y_s) = (n+1-\sum_{j=1}^s d_j,s)$, and $J(F)$ is the
Jacobian ideal.  Note the similarity with the second part of
Theorem~6.3.  It should be possible to prove a version of
\thetag{6.10} for complete intersections of ample hypersurfaces in an 
arbitrary complete simplicial toric variety.

\head \S7. Secondary fans and polytopes\endhead

	The {\it secondary polytope\/} was discovered in the year 1988 by
Gelfand, Kapranov and Zelevinsky \cite{GKZ1}, with further developments
by Billera, Filliman and Sturmfels \cite{BFS} and Oda and Park
\cite{OP}.  We will discuss the secondary polytope and its normal
fan, which is called the {\it secondary polytope\/} or {\it GKZ
decomposition\/}.  At the end of the section we will also mention the
{\it Gr\"obner fan\/} of a toric ideal.

	\subhead The secondary fan\endsubhead We begin with the
secondary fan, following \cite{OP}.  As noted in \S3, different toric
varieties can have closely related K\"ahler cones.  This happens
because different fans in $N_{\Bbb R}$ can share the same
$1$-dimensional cones.  To study the general case, fix a finite set of
primitive vectors $\Cal{B} \subset N$, where as usual $N\simeq{\Bbb
Z}^n$ and $M$ is its dual, and assume that $\roman{Cone}(\Cal{B}) =
N_{\Bbb R}$.  We will consider {\it all\/} complete fans $\Sigma$ in
$N_{\Bbb R}$ whose $1$-dimensional cone generators $\Sigma(1)$ lie in
$\Cal{B}$.

	Given $\Cal{B}$, we define $A_{\Cal{B}} \simeq {\Bbb
R}^{|\Cal{B}|-n}$ by the exact sequence
$$
0 \longrightarrow M_{\Bbb R} \overset\alpha\to\longrightarrow
{\textstyle{\bigoplus_{\rho \in \Cal{B}}}}{\Bbb R}\cdot e_\rho
\overset\beta\to\longrightarrow A_{\Cal{B}} \longrightarrow 0,\tag7.1
$$
where $\alpha(m) = \sum_{\rho\in\Cal{B}} \langle m,\rho\rangle
e_\rho$.  Also let $\Cal{B}^0 = \{\beta(e_\rho) : \rho \in
\Cal{B}\} \subset A_{\Cal{B}}$.  Then the pair
$(\Cal{B}^0,A_{\Cal{B}})$ is the {\it linear Gale transform\/} of
$(\Cal{B}, N_{\Bbb R})$. 
The cone
$$
A^+_{\Cal{B}} = \roman{Cone}(\Cal{B}^0) = \big\{{\textstyle{\sum_{\rho
\in \Cal{B}}}} a_\rho \beta(e_\rho) : a_\rho \ge 0\big\}\subset
A_{\Cal{B}}, 
$$
is strongly convex since $\roman{Cone}(\Cal{B}) = N_{\Bbb R}$.  From
\thetag{2.2}, we see that if $\Sigma$ is a fan with $\Sigma(1) =
\Cal{B}$, then $A_{\Cal{B}} =  A_{n-1}(X_\Sigma)\otimes{\Bbb R}$ and
$A^+_{\Cal{B}}$ is the cone generated by effective divisors on
$X_\Sigma$.

	Now suppose that $\Sigma$ is a simplicial projective fan with
$\Sigma(1) = \Cal{B}$.  From \S3, we have the cone
$\roman{cpl}(\Sigma) \subset A_{\Cal{B}}$ of convex classes, and it is
easy to see that $\roman{cpl}(\Sigma) \subset A^+_{\Cal{B}}$.
Furthermore, the interior of $\roman{cpl}(\Sigma)$ is the K\"ahler
cone of $X_\Sigma$.  In the example from \S3, these cones filled up
$A^+_{\Cal{B}}$ as we varied $\Sigma$.  In general, this doesn't
always happen, which is why we must allow fans
$\Sigma$ where $\Sigma(1)$ is a proper subset of $\Cal{B}$.  Here,
$A_{\Cal{B}}$ no longer equals $A_{n-1}(X_\Sigma)\otimes{\Bbb R}$, so
that the definition of convex needs to be modified.  We say that
$\sum_{\rho\in\Cal{B}} a_\rho \beta(e_\rho)$ is {\it convex\/} if for
every $\sigma \in \Sigma$, there is $m_\sigma \in M_{\Bbb R}$ such
that $\langle m_\sigma,\rho\rangle \ge -a_\rho$ for all $\rho \in
\Cal{B}$, with equality holding for $\rho \in \sigma$.  We again get a
strongly convex polyhedral cone $\roman{cpl}(\Sigma) \subset
A^+_{\Cal{B}}$ which has nonempty interior when $\Sigma$ is simplicial
and projective.  According to \cite{OP{\rm, Cor.~3.6}}, the cones
$\roman{cpl}(\Sigma)$ fit together as follows.

\proclaim{Theorem 7.1} If $\roman{Cone}(\Cal{B}) = N_{\Bbb R}$, then,
as $\Sigma$ ranges over all simplicial projective fans in $N_{\Bbb R}$
with $\Sigma(1) \subset \Cal{B}$, the cones $\roman{cpl}(\Sigma)$ and
their faces form a fan in $A_{\Cal{B}}$ whose support is
$A^+_{\Cal{B}}$.
\endproclaim
	
	The fan of this theorem is the {\it secondary fan\/} or {\it
GKZ decomposition\/} of $\Cal{B}$.  An example with two maximal cones
was given in \S3.  In general, the structure of the secondary fan is
quite interesting.  For example, two cones $\roman{cpl}(\Sigma)$ and
$\roman{cpl}(\Sigma')$ with a common codimension $1$ face are related
by a ``flop'' (as in the example from \S3) or by adding or subtracting
a single $1$-dimensional cone to $\Sigma$.  Also, faces on the
boundary of $A^+_{\Bbb R}$ correspond to certain ``degenerate'' fans.
This is all explained in \cite{OP{\rm, Sect.~3}}.

	The secondary fan has a nice relation to the moment map from
\S4.  If you look back at its construction, you'll see that the group
$G$ and the moment map $\mu_\Sigma$ of $X_\Sigma$ depend only on
$\Sigma(1)$.  Thus, given $\Cal{B} \subset N$ as above, we get the
moment map
$$
\mu_{\Cal{B}} : {\Bbb C}^{\Cal{B}} \longrightarrow A^+_{\Cal{B}}.
$$
By Theorem 4.1, if $\bold{a}$ is in the interior of
$\roman{cpl}(\Sigma)$, then
$$
\mu_{\Cal{B}}^{-1}(\bold{a})/G_{\Bbb R} \simeq X_\Sigma
$$
when $\Sigma$ is projective, simplicial and satisifies $\Sigma(1) =
\Cal{B}$.  In \cite{CK}, it is shown that this holds under the weaker
hypothesis $\Sigma(1) \subset \Cal{B}$ (with $\Sigma$ still projective
and simplicial).  Hence, the secondary fan gives a complete picture of
the toric varieties we can build via symplectic reduction from a given
moment map.  

	\subhead The secondary polytope\endsubhead Given a finite
subset $\Cal{A} \subset N_{\Bbb R} \simeq {\Bbb R}^n$, we can describe
a very interesting polytope using certain triangulations of its convex
hull $\Delta^* = \roman{Conv}(\Cal{A})$.  We will assume that
$\Delta^*$ is $n$-dimensional.  Then a {\it triangulation\/} of
$\Cal{A}$ is a triangulation $T$ of $\Delta^*$ such that the vertices
of each simplex in $T$ lie in $\Cal{A}$ (though not every element of
$\Cal{A}$ need be used).  Furthermore, the triangulation $T$ is {\it
regular\/} or {\it coherent\/} if there is a function $\psi : \Delta^*
\to {\Bbb R}$ which is affine on each simplex of $T$ and strictly
convex.  See \cite{GKZ1} or \cite{BFS} for precise definitions and
proofs of the results stated below.

	 If $T$ is a triangulation of $\Cal{A}$, we get the point
$$
\phi_T = \sum_{\sigma\in T}
\sum_{u\in\sigma}\roman{Vol}(\sigma)\,e_u \in
{\textstyle{\bigoplus_{u\in\Cal{A}}}} {\Bbb R}\cdot e_u.\tag7.2
$$
The convex hull of these points is the {\it secondary polytope\/} of
$\Cal{A}$, denoted $\Sigma(\Cal{A})$.  One can show that the dimension
of $\Sigma(\Cal{A})$ is $|\Cal{A}| -n-1$ and that its vertices are
precisely the points $\phi_T$ for which $T$ is a regular triangulation
of $\Cal{A}$.  This polytope has some very interesting combinatorial
properties \cite{GKZ1}, \cite{Stu2}.

	The secondary polytope also has a normal fan
$\Cal{N}(\Sigma(\Cal{A}))$, which is described using an affine version
of the Gale transform.  From $\Cal{A}$, we get the set
$\Cal{A}\times\{1\} \subset N_{\Bbb R}\otimes{\Bbb R}$, and as in
\thetag{7.1}, we can construct an exact sequence
$$
0 \longrightarrow M_{\Bbb R}\oplus{\Bbb R}
\overset\alpha\to\longrightarrow {\bigoplus_{u \in
\Cal{A}}} {\Bbb R}\cdot e_u 
\overset\beta\to\longrightarrow A_{\Cal{A}} \longrightarrow 0,\tag7.3 
$$
where $\alpha(m,\lambda) = \sum_{u\in \Cal{A}} (\langle m,u\rangle +
\lambda)e_u$.  Also set $\Cal{A}^0 = \{\beta(e_u) : u \in \Cal{A}\}$.
Then the pair $(\Cal{A}^0,A_{\Cal{A}})$ is the {\it affine Gale
transform\/} of $(\Cal{A},N_{\Bbb R})$.  Note that $\sum_{u\in\Cal{A}}
\beta(e_u) = 0$ and that $A_{\Cal{A}}$ has dimension $|\Cal{A}|-n-1$.

	Using the dual map $\beta^* : A_{\Cal{A}}^* \to \bigoplus_{u
\in \Cal{A}} {\Bbb R}\cdot e_u$, one can show that the image of
$\beta^*$ is parallel to the affine span of the secondary polytope
$\Sigma(\Cal{A})$.  It follows that the normal fan of
$\Sigma(\Cal{A})$ lives naturally in $A_{\Cal{A}}$.  This is the {\it
secondary fan of $\Cal{A}$\/}, denoted $\Cal{N}(\Sigma(\Cal{A}))$.
The maximal cones of the secondary fan correspond to vertices of
$\Sigma(\Cal{A})$ and hence to regular triangulations of $\Cal{A}$.
Note also that $\Cal{N}(\Sigma(\Cal{A}))$ is a complete fan in
$A_{\Cal{A}}$ since it comes from a polytope.

	\subhead The enlarged secondary fan\endsubhead The secondary
fan $\Cal{N}(\Sigma(\Cal{A}))$ might seem rather different from the
secondary fan defined earlier.  To see the relation, let $\Cal{B} =
\Sigma(1)$, where $\Sigma$ is a projective fan, not necessarily
simplicial.  Then, as before, we get the secondary fan of $\Cal{B}$,
whose support is the strongly convex cone $A^+_{\Cal{B}} \subset
A_{\Cal{B}}$.  Following \cite{AGM2}, we can enlarge this to a
complete fan as follows.  Fix an ample divisor $D = \sum_\rho a_\rho
D_\rho$ with $a_\rho > 0$ for all $\rho$.  This implies that $0$ is an
interior point of the corresponding polytope $\Delta \subset M_{\Bbb
R}$.  Then the {\it dual polytope\/} $\Delta^\circ \subset N_{\Bbb R}$
is defined by
$$
\Delta^\circ = \{u \in N_{\Bbb R} : \langle m,u\rangle \ge -1\ \hbox{for
all}\ m \in \Delta\}.
$$
Note that $0$ is an interior point of $\Delta^\circ$, though
$\Delta^\circ$ need not be integral.  Also, the normal fan $\Sigma$ of
$\Delta$ is obtained by taking cones over proper faces of
$\Delta^\circ$.  In particular, the vertices of $\Delta^\circ$ are
$(1/a_\rho)\rho$ for $\rho
\in \Cal{B} = \Sigma(1)$.  Now let 
$$
\Cal{A} = \roman{Vert}(\Delta^\circ)\cup\{0\} = \{(1/a_\rho)\rho :
\rho \in \Sigma(1)\} \cup\{0\}.
$$ 
We leave it to the reader to show that there is a natural isomorphism
$A_{\Cal{A}} \simeq A_{\Cal{B}}$ which carries $\Cal{A}^0$ to
$\Cal{B}^0\cup\{-\sum_{\rho\in\Cal{B}} a_\rho \beta(e_\rho)\}$.  Under
this isomorphism, we can compare the secondary fans of $\Cal{A}$ and
$\Cal{B}$ as follows.

\proclaim{Proposition 7.2} Let $\Cal{A}$ and $\Cal{B}$ be as above.
Then: 
\roster
\item Under the isomorphism $A_{\Cal{A}} \simeq A_{\Cal{B}}$, the cone
of $\Cal{N}(\Sigma(\Cal{A}))$ given by a regular triangulation
$T$ of $\Cal{A}$ lies in $A^+_{\Cal{B}}$ if and only if $0$ is
contained in every maximal simplex of $T$.
\item For such a triangulation $T$, let $\Sigma'$ is the fan obtained
by taking cones over simplices of $T$.  Then $\Sigma'$ is a projective
simplicial fan with $\Sigma'(1) \subset \Cal{B}$, and the cone
$\roman{cpl}(\Sigma')$ is the cone of $\Cal{N}(\Sigma(\Cal{A}))$
given by to $T$.
\item Conversely, every maximal cone $\roman{cpl}(\Sigma') \subset
A^+_{\Cal{B}}$ comes from such a triangulation $T$.  We get $T$ by
intersecting $\Delta^\circ$ with the cones in $\Delta^\circ$.
\endroster
\endproclaim

	In this situation, we call $\Cal{N}(\Sigma(\Cal{A}))$ the {\it
enlarged secondary fan\/} of $\Cal{B}$.  A consequence of this
proposition is that there is a bijective correspondence between
maximal cones of the secondary fan of $\Cal{B}$ and regular
triangulations of $\Cal{A}$ whose maximal cones all contain $0$.

\demo{Example} The toric variety $X = X_\Sigma$ from \S3 has cone
generators 
$$
e_0 = (0,0,-2),\ e_1 = (1,1,1),\ e_2 = (1,-1,1),\ e_3 = (-1,-1,1),\
e_4 = (-1,1,1)
$$
and maximal cones $\sigma_{1234}$, $\sigma_{012}$, $\sigma_{023}$,
$\sigma_{034}$ and $\sigma_{041}$.  Let $\Delta$ be polytope
associated to the anticanonical divisor $D_0 + D_1 + D_2 + D_3 + D_4$,
which is ample.  One can compute that the dual polytope of
$\Delta$ is
$$
\Delta^\circ = \roman{Conv}(e_0,e_1,e_2,e_3,e_4).\tag7.4
$$
In \S8, we will see that $\Delta$ and $\Delta^\circ$ are examples of {\it
reflexive polytopes\/}.

	Using $\Cal{A} = \roman{Vert}(\Delta^\circ) \cup\{0\} =
\{e_0,e_1,e_2,e_3,e_4,0\}$, one sees that there are exactly
four triangulations of $\Cal{A}$, which means that the enlarged
secondary fan has four maximal cones.  Two of the triangulations have
$0$ as a vertex, and the corresponding cones were described in \S3.
We leave it to the reader to determine the other two cones in the fan
and the corresponding triangulations of $\Cal{A}$.
\enddemo

	One can also describe the enlarged secondary fan using the
total space of the line bundle over $X_\Delta$ given by $-\sum_\rho
a_\rho D_\rho$ (see \cite{AGM2}).  This is closely related to an
alternate approach to the secondary fan $\Cal{N}(\Sigma(\Cal{A}))$
which appears in \cite{OP}.  We will see in \S10 that the
enlarged secondary fan and the secondary polytope are used in mirror
symmetry.  

	Other references for secondary polytopes and the
GKZ-decomposition are \cite{Loe} and \cite{Par}.  A significant
generalization of the secondary polytope, called the {\it fiber
polytope\/}, is described in \cite{BS} and \cite{Zie}.
	
	\subhead The Gr\"obner fan\endsubhead If our finite set
$\Cal{A}$ is integral, i.e., $\Cal{A} = \{u_1,\ldots,u_\ell\} \subset
N$, then we can use term orders on a toric ideal to create a fan
closely related to the secondary fan of $\Cal{A}$.  The basic idea is
that $\Cal{A}\times\{1\} \subset N\oplus{\Bbb Z}$ determines the toric
ideal $I_{\Cal{A}\times\{1\}} \subset {\Bbb C}[x_1,\dots,x_\ell]$.
This ideal defines the affine toric variety
${}^a\!X_{\Cal{A}\times\{1\}}$, which by \S5 is the affine cone over
$X_{\Cal{A}} \subset {\Bbb P}^{\ell-1}$.  For this reason, we write
the toric ideal as $I_{\Cal{A}}$ instead of $I_{\Cal{A}\times\{1\}}$.
By Lemma 5.5, we can express this ideal as
$$
I_{\Cal{A}} = \langle \bold{x}^{\bold{a}^+} -
\bold{x}^{\bold{a}^-} : \bold{a} = (a_1,\dots,a_\ell) \in {\Bbb
Z}^\ell,\ {\textstyle{\sum_{i=1}^\ell}} a_im_i = 0,\
{\textstyle{\sum_{i=1}^\ell}} a_i= 0\rangle.\tag7.5
$$

	We can use elements $\omega \in {\Bbb R}^\ell$ to create
initial ideals of $I_{\Cal{A}}$ (in the sense of Gr\"obner theory) as
follows: if $f = \sum_{\bold{a}} c_{\bold{a}} \bold{x}^{\bold{a}}$,
define ${\hbox{\smc lt}}_\omega(f) = \sum_{\omega\cdot\bold{a}\
\roman{maximal}} c_{\bold{a}} \bold{x}^{\bold{a}}$, and then set
$$
{\hbox{\smc lt}}_\omega(I_{\Cal{A}}) = \langle {\hbox{\smc
lt}}_\omega(f) : f \in I_{\Cal{A}}\rangle.
$$
For generic $\omega$, this is the initial ideal of $I_{\Cal{A}}$ for
some term order, and one can show that every initial ideal arises in
this way \cite{Stu2}.  The interesting aspect is that different
$\omega$'s can give the same initial ideal.  We define
$$
\omega \sim \omega' \iff {\hbox{\smc lt}}_\omega(I_{\Cal{A}}) = {\hbox{\smc
lt}}_{\omega'}(I_{\Cal{A}}).
$$
One can prove that each equivalence class is a relatively open convex
polyhedral cone.  Furthermore, if $\alpha$ is the map from
\thetag{7.3}, then $\omega + \alpha(m,\lambda) \sim \omega$ for all
$(m,\lambda) \in M_{\Bbb R}\oplus{\Bbb R}$.  Hence it makes sense to
take the quotient by the image of $\alpha$.  This leads to the
following result.

\proclaim{Proposition 7.3} Under the map $\beta : {\Bbb R}^\ell \to
A_{\Cal{A}}$ from \thetag{7.3}, we have:
\roster
\item The closures of the images of the equivalence classes form a
complete fan in $A_{\Cal{A}}$ called the {\rm Gr\"obner fan}.
\item The Gr\"obner fan of $\Cal{A}$ refines the secondary fan
$\Cal{N}(\Sigma(\Cal{A}))$. 
\endroster
\endproclaim

	This is proved in \cite{Stu2}.  One can also show that the
Gr\"obner fan is the normal fan of a polytope called the {\it state
polytope\/}.  Algorithms for computing the state polytope and
Gr\"obner fan are given in \cite{Stu2}, and the reader may also wish
to consult \cite{MR}.  Although the Gr\"obner fan of $I_{\Cal{A}}$ may
be strictly finer than the secondary fan of $\Cal{A}$, there are
situations where the Gr\"obner fan is very useful.  An example from
mirror symmetry can be found in \cite{HLY1}.

\head \S8. Reflexive polytopes and Calabi-Yau hypersurfaces\endhead

	In this section, we will create some interesting families of
Calabi-Yau varieties using toric geometry.  The key idea will be
Batyrev's notion of {\it reflexive polytopes\/}.  At the end of the
section we will also consider complete intersections and
nef-partitions.

	\subhead Singular Calabi-Yau varieties\endsubhead The quintic
threefold in ${\Bbb P}^4$ is one of the best-known examples of a
Calabi-Yau manifold.  However, there are many contexts (including
toric geometry) where {\it singular\/} Calabi-Yau varieties arise
naturally.  Hence we define Calabi-Yau as follows.  Let $Y$ be a
$d$-dimensional normal projective complex variety.  Assume 
$\Cal{O}_Y(K_Y) = \widehat\Omega_Y^d \simeq \Cal{O}_Y$ and 
$$
H^1(Y,\Cal{O}_Y) \simeq \cdots \simeq H^{d-1}(Y,\Cal{O}_Y) \simeq
\{0\}. 
$$
Then $Y$ is {\it canonical Calabi-Yau\/} if it has at worst canonical
singularities.  Furthermore, $Y$ is {\it minimal Calabi-Yau\/} if it
has at worst ${\Bbb Q}$-factorial terminal singularities.

	Since canonical singularities are Cohen-Macaulay and a
Calabi-Yau has trivial canonical bundle, a singular Calabi-Yau is
always Gorenstein.  Thus its singularities are either Gorenstein
canonical or ${\Bbb Q}$-factorial Gorenstein terminal.

	To understand why these singularities are appropriate for
Calabi-Yau varieties, recall that a singularity is {\it canonical\/}
if there is a local resolution of singularities $f : \widehat Y \to Y$
such that
$$
K_{\widehat Y} = f^*(K_Y) + \sum_i a_i E_i,
$$
where the sum is over the exceptional divisors $E_i$ of $f$ and $a_i
\ge 0$.  If in addition we have $a_i > 0$ for all $i$, then the
singularity is {\it terminal\/}.  Thus having terminal singularities
means that any resolution must change the canonical class.  In
particular, if $Y$ is a singular minimial Calabi-Yau, then its
resolutions cease to be Calabi-Yau since the canonical class is no
longer trivial.

	As for canonical singularities, there are many situations
(including threefolds and the Calabi-Yaus to be constructed below)
where having canonical singularities implies the existence of a
partial resolution $f: \widehat Y \to Y$ such that $K_{\widehat Y} =
f^*(K_Y)$ (so the canonical class doesn't change) and $\widehat Y$ has
${\Bbb Q}$-factorial terminal singularities (so any further resolution
changes the canonical class).  Thus, if we start from a singular
canonical Calabi-Yau $Y$, then $\widehat Y$ (if it exists) is a
minimal Calabi-Yau which is as close as possible to being smooth while
remaining Calabi-Yau.  In Batyrev's terminology, $\widehat Y$ is
called a {\it maximal projective crepant partial desingularization\/}
of $Y$ (a MPCP-desingularization for short).

	For more background on canonical and terminal singularities
(including some nice examples), see \cite{Rei5}.

	\subhead Reflexive polytopes and Fano toric
varieties\endsubhead Following Batyrev \cite{Bat1}, we say that a
$n$-dimensional integral convex polytope $\Delta \subset M_{\Bbb R}$
is {\it reflexive\/} if it contains the origin as an interior point
and if its dual polytope
$$
\Delta^\circ = \{u \in N_{\Bbb R} : \langle m,u\rangle \ge -1\
\hbox{for all}\ m \in \Delta\} \subset N_{\Bbb R}
$$
is also integral.  Since $(\Delta^\circ)^\circ = \Delta$, reflexive
polytopes always come in pairs.  It also follows that each facet $F$
of a reflexive polytope $\Delta$ is defined by an equation $\langle
u_F,m\rangle = -1$ for some $u_F \in N$.  This easily implies that $0$
is the only point of $M$ in the interior of $\Delta$.  The polytopes
$\Delta$ and $\Delta^\circ$ from \thetag{7.4} are an example of a
reflexive polytope and its dual.

	It is known (see \cite{Bat1}) that in each dimension, there
are only finitely many reflexive polytopes up to unimodular
equivalence, and work is underway to classify {\it all\/} reflexive
polytopes in dimension $4$ (see \cite{Ska}).  

	Given a reflexive polytope $\Delta$, we get a toric variety
$X_\Delta$.  Since the facets of $\Delta$ are given by $\langle
u_F,m\rangle = -1$, it follows easily that the divisor on $X_\Delta$
determined by $\Delta$ is precisely the anticanonical divisor $-K =
\sum_\rho D_\rho$.  Also, the anticanonical divisor is ample, which
tells us that $X_\Delta$ is a Fano variety, and it is Gorenstein
since $K$ is obviously Cartier.  Conversely, one can show that every
Gorenstein Fano toric variety arises from a reflexive polytope (see
\cite{Bat1}).

	\subhead Calabi-Yau hypersurfaces\endsubhead Now suppose that
$Y \subset X_\Delta$ is a general anticanonical hypersurface in a
Gorenstein Fano toric variety coming from a reflexive polytope
$\Delta$.  Then the adjunction formula
$$
\widehat\Omega_Y^{n-1} \simeq
\widehat\Omega_{X_\Delta}^n(-K)\otimes\Cal{O}_Y 
$$
shows that $\widehat\Omega^{n-1}_Y \simeq \Cal{O}_Y$.  Since
$\Cal{O}_{X_\Delta}(-Y) \simeq \Cal{O}_{X_\Delta}(K) = 
\widehat\Omega_{X_\Delta}^n$, we also get an exact sequence
$$
0 \longrightarrow \widehat\Omega_{X_\Delta}^n \longrightarrow
\Cal{O}_X \longrightarrow \Cal{O}_Y \longrightarrow 0,
$$
which makes it easy to show that $H^i(Y,\Cal{O}_Y) = 0$ for $0 < i <
n-1$ (this uses Serre-Grothendieck duality and the vanishing of
$H^i(X_\Delta,\Cal{O}_{X_\Delta})$ for $i > 0$).  Finally, Bertini
theorems show that $Y$ has at most Gorenstein toric singularities,
which are known to be canonical.  It follows that $Y$ is a canonical
Calabi-Yau variety.

	When $Y$ is singular, we would like to desingularize it as
much as possible while remaining Calabi-Yau.  As one might expect,
toric geometry tells us exactly what to do.  We've seen that the cone
generators of the normal fan are the vertices of $\Delta^\circ$.  We will
consider projective simiplicial fans $\Sigma$ which refine the normal
fan of $\Delta$ and whose cone generators satisfy
$$
\Sigma(1) = N\cap \Delta^\circ - \{0\}.\tag8.1
$$
To see why this condition is relevant, remember that one
desingularizes $X_\Delta$ by subdividing cones into simplicial ones
and then adding new cone generators until we get smooth cones.  For
example, adding a new cone generator $\rho$ gives a birational map of
toric varieties $f : X \to X_\Delta$, and if $\rho$ lies in the cone
over the facet $\langle u,m_F\rangle = -1$ of $\Delta^\circ$, then
using the techniques of \cite{Rei5}, one obtains
$$
K_{X} = f^*(K_{X_\Delta}) - (\langle \rho,m_F\rangle +
1)D,\tag8.2
$$
where $D$ is the exceptional divisior (and is the divisor of $X$
corresponding to $\rho$).  Thus, as long as we add new cone
generators coming from $N\cap\Delta^\circ$, we don't change the canonical
class.  But once we've used up all of $N\cap \Delta^\circ -
\{0\}$, as $\Sigma$ does, then any further $\rho$'s must lie
outside of $\Delta^\circ$ (since reflexive implies $0$ is the only integer
interior point).  Hence $\langle\rho,m_F\rangle < -1$, and from
\thetag{8.2}, it follows that we have terminal singularities.  Since
$\Sigma$ is simplicial, we see that $X_\Sigma$ is ${\Bbb
Q}$-factorial.  Thus we have the following theorem (see \cite{Bat1}
for details).

\proclaim{Theorem 8.1} Let $\Delta \subset M_{\Bbb R}$ be a reflexive
polytope, and let $\Sigma$ be a fan in $N_{\Bbb R}$ which refines the
normal fan of $\Delta$ and satisfies \thetag{8.1}.  Then:
\roster
\item The general anticanonical hypersurface $Y \subset X_\Delta$ is a
canonical Calabi-Yau variety.
\item The general anticanonical hypersurface $\widehat Y \subset
X_\Sigma$ is a minimal Calabi-Yau variety.  Furthermore, $\widehat Y$
is a MPCP-desingularization of its image $Y \subset X_\Delta$ under
the map $X_\Sigma \to X_\Delta$.
\endroster
\endproclaim

	When $\Delta$ is a $4$-dimensional reflexive polytope, one can
show that the MPCP-desingularization $\widehat Y \subset X_\Sigma$ is
smooth Calabi-Yau threefold (see \cite{Bat1}).

	It follows that for each fan $\Sigma$ as in the statement of
the theorem, we get a family of minimal Calabi-Yau varieties.  In
terms of what we considered in \S7, these fans correspond to certain
maximal cones in the secondary fan of the set $\Cal{B} =
N\cap\Delta^\circ-\{0\} \subset N$.  

	Of course, since $\Delta$ gives the family $\widehat Y
\subset X_\Sigma$ of minimal Calabi-Yaus, we can use the reflexive
polytope $\Delta^\circ$ to construct a ``dual'' family $\widehat
Y^\circ \subset X_{\Sigma^\circ}$ of minimal Calabi-Yau varieties,
where $\Sigma^\circ(1) = M\cap\Delta-\{0\}$.  We will see in \S10 that
$\widehat Y^\circ \subset X_{\Sigma^\circ}$ is conjectured to be the
mirror family of $\widehat Y \subset X_\Sigma$.  Because of the
previous paragraph, the situation is complicated by the multiple
choices for $\Sigma$ and $\Sigma^\circ$.  This is related to the idea
of {\it multiple mirrors\/}, also to be studied in \S10.

	There are also some nice formulas for Hodge numbers.  The
MPCP-desingular\-iza\-tion $\widehat Y \subset X_\Sigma$ is
simplicial, so that $H^*(\widehat Y)$ has a pure Hodge structure.  But
we can't apply the results of \S6 directly since $\widehat Y$ may fail
to be ample (this happens, for example, when $X_\Delta = {\Bbb
P}(1,1,2,2,2)$).  However, we have the following results from
\cite{Bat1}.

\proclaim{Theorem 8.2} Let $\Delta$ be a reflexive polytope of
dimension $n \ge 4$ and let $\Sigma$ be a fan as in Theorem 8.1.  If
$\widehat Y \subset X_\Sigma$ is a general anticanonical hypersurface,
then
$$
\align
h^{n-2,1}(\widehat Y) &= l(\Delta) - n - 1 -
   \sum_{\roman{codim}\,\Gamma = 1} l^*(\Gamma) +  
   \sum_{\roman{codim}\,\Gamma = 2} l^*(\Gamma) l^*(\Gamma^\circ),\\
\intertext{where $\Gamma$ is a face of $\Delta$, $\Gamma^\circ$ is the
corresponding dual face of $\Delta^\circ$ and, as in \S6, $l(\Gamma)$
(resp.\ $l^*(\Gamma)$) is the number of integer (resp.\ interior
integer) points in $\Gamma$.  Also,}
h^{1,1}(\widehat Y) &= l(\Delta^\circ) - n - 1 -
\sum_{\roman{codim}\,\Gamma^\circ = 1} l^*(\Gamma^\circ) +  
\sum_{\roman{codim}\,\Gamma^\circ = 2} l^*(\Gamma^\circ) l^*(\Gamma),
\endalign
$$
where now $\Gamma^\circ$ is a face of $\Delta^\circ$ and $\Gamma$ is
the dual face.
\endproclaim

	The remarkable symmetry evident between the formulas for
$h^{1,1}$ and $h^{n-2,1}$ will be important when we discuss mirror
symmetry in \S10.  A different approach to the study of $h^{1,1}$ can
be found in \cite{Roa1}.  

	Although the above construction seems completely natural, it
can take some thought to find the reflexive polytope.  For example, a
degree $7$ hypersurface in ${\Bbb P}(1,1,1,2,2)$ is Calabi-Yau, yet
the simplex giving ${\Bbb P}(1,1,1,2,2)$ in not reflexive, nor is the
hypersurface Cartier.  Here, the reflexive polytope $\Delta$ is (up to
translation) the Newton polytope of all monomials of degree $7$, and
the toric variety $X_\Delta$ is a blow-up of ${\Bbb P}(1,1,1,2,2)$.
Some early evidence for mirror symmetry came from a list of 7555
Calabi-Yau hypersurfaces in weighted projective spaces.  The list
exhibited an incomplete duality which was only fully understood after
the definition of reflexive polytope (see
\cite{CdK}).

	\subhead Calabi-Yau complete intersections\endsubhead We next
generalize the above construction to create families of Calabi-Yau
complete intersections.  In \S6, we discussed some basic facts about
complete intersections in toric varieties, though the special features
of the Calabi-Yau case are due to \cite{LBor}, with subsequent work by
\cite{BB1--3}.  Details of what follows can be found in these
references.

	The basic way to describe a Calabi-Yau complete intersection
in a toric variety is through the idea of a nef-partition.  Suppose
that $\Delta$ is a $n$-dimensional reflexive polytope, so that
$\roman{Vert}(\Delta^\circ)$ is the set of cone generators of the fan of
$X = X_\Delta$.  Then a {\it nef-partition\/} is a disjoint union
$$
\roman{Vert}(\Delta^\circ) = E_1\cup\cdots\cup E_r\tag8.3
$$
such that the divisors $D_j = \sum_{\rho\in E_j} D_\rho$ are Cartier
and nef (i.e., generated by global sections).  Equivalently, a
nef-partition is a partition of the anticanonical divisor $-K =
\sum_\rho D_\rho$ into a sum of nef Cartier divisors.  Furthermore, if
$\Delta_j$ is the polytope associated to the divisor $D_j$, then it
follows that
$$
\Delta = \Delta_1 + \cdots + \Delta_r.\tag8.4
$$
We will always assume that $n-r \ge 1$.

	In this situation, let $Y_j \in |D_j|$ be generic.  Then the
complete intersection
$$
V = Y_1\cap\cdots\cap Y_r \subset X_\Delta
$$
is a canonical Calabi-Yau.  Furthermore, as in Theorem 8.1, suppose
that $\Sigma$ is a fan refining the normal fan of $\Delta$ with the
property that $\Sigma(1) = N\cap\Delta^\circ-\{0\}$.  We can regard
$D_j$ as divisors on $X_\Sigma$, and if $\widehat Y_j$ is a general
member of $|D_j|$, then, as shown in \cite{BB3}, the complete
intersection
$$
\widehat V = \widehat Y_1\cap\cdots\cap \widehat Y_r \subset
X_\Sigma\tag8.5 
$$
is a minimal Calabi-Yau variety which is a MPCP-desingularization of
the corresponding $V \subset X_\Delta$.  There are formulas similar to
those of Theorem 8.2 for the Hodge numbers $h^{n-r-1,1}(\widehat V)$
and $h^{1,1}(\widehat V)$.

	In \S10, we will see that the nef-partition giving the family
of Calabi-Yau complete intersections $\widehat V \subset X_\Sigma$
naturally determines a ``dual'' family of Calabi-Yau complete
intersections.  This construction is used in mirror symmetry.

\head \S9. Resultants, discriminants and hypergeometric functions\endhead

	This section will discuss briefly those aspects of resultants,
discriminants and hypergeometric functions which pertain to toric
varieties.  Basic references are the book \cite{GKZ1} and, for
hypergeometric functions, the papers \cite{GKZ2} and \cite{GKZ3}.

	\subhead $\Cal{A}$-Resultants\endsubhead Given a finite set of
exponents $\Cal{A} \subset {\Bbb Z}^n$, we get the vector space
$L(\Cal{A})$ of Laurent polynomials
$$
f = \sum_{m \in \Cal{A}} c_m \bold{t}^m.
$$
For simplicity, we will assume that $\Cal{A}$ is affinely independent
in ${\Bbb Z}^n$.  Then let $\nabla_{\Cal{A}} \subset L(\Cal{A})^{n+1}$
be the Zariski closure of those $(n+1)$-tuples of polynomials
$(f_0,\dots,f_n) \in L(\Cal{A})^{n+1}$ where there is $\bold{t} \in
({\Bbb C}^*)^n$ such that $f_0(\bold{t}) = \cdots = f_n(\bold{t}) =
0$.  One can show that $\nabla_{\Cal{A}}$ is a hypersurface in
$L(\Cal{A})^{n+1}$ and hence is defined by a polynomial $R_{\Cal{A}} =
0$.  Thus $R_{\Cal{A}}$ is a polynomial in the coefficients $c_{i,m}$
of $f_i$, and in fact $R_{\Cal{A}} \in {\Bbb Z}[c_{i,m}]$.  This
polynomial is the {\it $\Cal{A}$-resultant\/}.

	For example, if $\Cal{A} = \{0,1,\ldots,d\} \subset {\Bbb Z}$,
then $R_{\Cal{A}}$ is the usual resultant of two polynomials $f_0,f_1$
of degree $d$ in one variable.  Many other examples can be found in
\cite{GKZ3} and \cite{Stu4}.  In general, the $\Cal{A}$-resultant
is sometimes called a {\it sparse resultant\/} since only the
exponents in $\Cal{A}$ are used.  It is also possible to allow the
polynomials $f_0,\dots,f_n$ to have different exponents---this leads
to what is known as the {\it
$(\Cal{A}_0,\dots,\Cal{A}_n)$-resultant\/} or {\it mixed sparse
resultant\/}.  Algorithms for computing sparse resultants are
described in \cite{CE} and \cite{Stu4}.

	\subhead Chow forms\endsubhead We can think of the
$\Cal{A}$-resultant in toric terms as follows.  If $\Cal{A}$ has
$\ell$ elements, then we get the (possibly non-normal) toric
variety $X_{\Cal{A}} \subset {\Bbb P}^{\ell-1}$ as in \thetag{5.6}.
This projective variety has a {\it Chow form\/} $R_{X_{\Cal{A}}}$,
which is a polynomial in the Pl\"ucker coordinates $[m_0,\dots,m_n]$
that vanishes precisely on those subspaces $L \in G(\ell,n-\ell-1)$
where ${\Bbb P}(L)\cap X_{\Cal{A}} \ne \emptyset$.  

	It is customary to rewrite the Chow form as follows.  If we
regard $L$ as being defined by $n+1$ linear forms, then the
coefficients of these forms give a $(n+1)\times\ell$ matrix
$(c_{i,m})$, where we use elements $m \in \Cal{A}$ to label the
coordinates of ${\Bbb P}^{\ell-1}$.  Then, replacing the Pl\"ucker
coordinate $[m_0,\dots,m_n]$ by the bracket polynomial
$[m_0,\dots,m_n] = \det(c_{i,m_j})$, we get a polynomial in the
$c_{i,m}$ which is easily seen to coincide with the
$\Cal{A}$-resultant $R_{\Cal{A}}$.

	\subhead Chow polytopes and secondary polytopes\endsubhead
Using the Chow form, we can create the {\it Chow polytope\/} as
follows.  If we express $R_{X_{\Cal{A}}}$ as a polynomial in the
brackets $[\sigma] = [m_0,\dots,m_n]$, then a term $T$ in
$R_{X_{\Cal{A}}}$ is of the form $T = c\Pi_\sigma
[\sigma]^{a_\sigma}$, $c \ne 0$, and we assign $T$ the weight
$$
\phi_T = \sum_\sigma \sum_{m \in [\sigma]} a_\sigma e_{m} \in
{\textstyle{\bigoplus_{m\in\Cal{A}}}} {\Bbb R}\cdot e_m. 
$$
The Chow polytope of $X_{\Cal{A}}$ is defined to be the convex hull of
these points $\phi_T$.

	The remarkable fact is that the Chow polytope of $X_{\Cal{A}}$
is {\it exactly\/} the secondary polyotope $\Sigma(\Cal{A})$ from \S7
(which is the convex hull of the points $\phi_T$ defined in
\thetag{7.2}).  This is proved in \cite{GKZ1{\rm, Ch.~9}}.  We know
that the vertices of $\Sigma(\Cal{A})$ correspond to regular
triangulations of $\Cal{A}$, and given such a triangulation, there is
a precise formula for the corresponding term of $R_{X_{\Cal{A}}}$.
One can also formulate this in terms of certain toric degenerations of
$X_{\Cal{A}}$ (see \cite{KSZ1}).  Some explicit examples can be found
in \cite{Stu4}, and \cite{Stu3} has similar results concerning the
$(\Cal{A}_0,\dots,\Cal{A}_n)$-resultant.

	\subhead $\Cal{A}$-Discriminants\endsubhead Given $X_{\Cal{A}}
\subset {\Bbb P}^{\ell-1}$ as above, let $\nabla_{\Cal{A}} \subset
L(\Cal{A})$ be the Zariski closure of those Laurent polynomials $f$
for which the affine hypersurface $Z_f \subset ({\Bbb C}^*)^n$ is
singular.  Then $\nabla_{\Cal{A}}$ is the affine cone over the dual
variety $X^\vee_{\Cal{A}}$.  If $\nabla_{\Cal{A}}$ has codimension 1,
we define the {\it $\Cal{A}$-discriminant\/} $\Delta_{\Cal{A}}$ to be
the defining equation of $\nabla_{\Cal{A}}$ (or, equivalently, of
$X^\vee_{\Cal{A}}$).  We set $\Delta_{\Cal{A}} = 1$ in all other
cases.  Note that $\Delta_{\Cal{A}}$ is a polynomial in the
coefficients $c_m$ of $f = \sum_{m\in\Cal{A}} c_m \bold{t}^m$, and,
in fact, $\Delta_{\Cal{A}} \in {\Bbb Z}[c_m]$.

\demo{Examples} 1. When $\Cal{A}$ consists of all non-negative
exponent vectors of degree at most $d$, $\Delta_{\Cal{A}}$ is the
usual discriminant of a homogeneous polynomial of degree $d$.  More
precisely, if $F$ is the homogenization of $f \in L(\Cal{A})$, then
$\Delta_{\Cal{A}}(f) = \roman{Disc}(F)$.

	2. Consider all linear combinations of the monomials
$1,x,\dots,x^p,y,yx,\dots,yx^q$.  Then $\Delta_{\Cal{A}}$ is the usual
resultant $R(f,g)$ of polynomials of degrees $p$ and $q$ respectively.
This is an example of the Cayley trick and can be used to express any
$(\Cal{A}_0,\dots,\Cal{A}_n)$-resultant as a discriminant (see
\cite{GKZ{\rm, Prop.~1.7 of Ch.~9}}).

	3. An example where $\Delta_{\Cal{A}} = 1$ is given by the
Segre embedding ${\Bbb P}^1\times{\Bbb P}^2 \to {\Bbb P}^3$.  More
generally, same is true for the Segre embedding of ${\Bbb
P}^n\times{\Bbb P}^m$ whenever $n \ne m$.
\enddemo

	\subhead Principal $\Cal{A}$-determinants\endsubhead Given a
Laurent polynomial $f
\in L(\Cal{A})$, the {\it principal $\Cal{A}$-determinant\/}
$E_{\Cal{A}}(f)$ is defined to be the resultant
$$
E_{\Cal{A}}(f) = R_{\Cal{A}}(f,t_1\partial f/\partial t_1,\dots,
t_n\partial f/\partial t_n).
$$
Note that the resultant makes sense since $t_i\partial f/\partial t_i
\in L(\Cal{A})$.

	The nicest fact about $E_{\Cal{A}}$ is that its Newton
polytope is {\it precisely\/} the secondary polytope $\Sigma(\Cal{A})$
(see \cite{GKZ{\rm, Thm.~1.4 of Ch.~10}}).  One way to understand this
result is to observe that $E_{\Cal{A}}$ is the polynomial obtained
from the Chow form $R_{X_{\Cal{A}}}$ under the specialization which
sends the bracket $[\sigma]$ to $\pm\roman{Vol}(\sigma)
\Pi_{m\in\sigma} c_m$ (where $\pm\roman{Vol}(\sigma)$ is the signed
volume of the simplex spanned by $m \in \sigma$).

	When $X_{\Cal{A}}$ is smooth, there is an especially nice
relation between $\Delta_{\Cal{A}}$ and $E_{\Cal{A}}$.

\proclaim{Theorem 9.1} If $X_{\Cal{A}} \subset {\Bbb P}^{\ell-1}$ is
smooth and $Q = \roman{Conv}(\Cal{A})$, then  
$$
E_{\Cal{A}}(f) = \pm{\textstyle{\prod_{\Gamma \subset Q}}}
\Delta_{\Cal{A}\cap\Gamma}(f_{|\Gamma}), 
$$
where $\Gamma \subset Q$ are the nonempty faces and
$f_{|\Gamma} = \sum_{m\in\Cal{A}\cap\Gamma} c_m \bold{t}^m$.
\endproclaim

	When $X_{\Cal{A}}$ is singular, there is a more complicated
but very elegant formula relating $\Delta_{\Cal{A}}$ and
$E_{\Cal{A}}$---see \cite{GKZ{\rm, Ch.~10}} for the details.

	Another way to see the relation between $\Delta_{\Cal{A}}$ and
$E_{\Cal{A}}$ is to study the projective hypersurface $Y_f \subset
X_{\Cal{A}}$ defined by $f \in L(\Cal{A})$.  Then we have the
following result.

\proclaim{Proposition 9.2} Given $X_{\Cal{A}} \subset {\Bbb
P}^{\ell-1}$, we have:
\roster
\item If $X_{\Cal{A}}$ is smooth and $\Delta_{\Cal{A}}$ is
nonconstant, then $Y_f \subset X_{\Cal{A}}$ is smooth if and only if
$\Delta_{\Cal{A}}(f) \ne 0$.
\item If $E_{\Cal{A}}$ is nonconstant, then $Y_f \subset
X_{\Cal{A}}$ is nondegenerate (as defined in \S6) if and only if
$E_{\Cal{A}}(f) \ne 0$.
\endroster
\endproclaim

	It is likely that the first part of this proposition remains
true in the simplicial case, so that $Y_f$ should be quasi-smooth in
the sense of \S6 if and only if $\Delta_{\Cal{A}}(f) \ne 0$.  However,
no proof has appeared in print.

	As a final comment, note that resultants can be computed in
terms of discriminants (using the Cayley trick) and vice versa (using
$E_{\Cal{A}}$ and Theorem 9.1 in the smooth case).  Historically,
$\Cal{A}$-discriminants were first discovered in the context of
$\Cal{A}$-hypergeometric functions and led to the definition of
$\Cal{A}$-resultants.  In practice, resultants are easier to compute
(see, for example, \cite{Stu4}).

	\subhead $\Cal{A}$-Hypergeometric functions\endsubhead To
define the $\Cal{A}$-hypergeometric equations, consider $\Cal{A} =
\{m_1,\dots,m_\ell\} \subset {\Bbb Z}^n$.  For simplicity, we will
assume that ${\Bbb Z}^n$ is affinely generated by $\Cal{A}$, which
means that ${\Bbb Z}^{n+1}$ is generated by $\Cal{A}\times\{1\}$.
Then, for each $\bold{a} = (a_1,\dots,a_n) \in {\Bbb Z}^\ell$
satisfying the conditions
$$
\sum_{i=1}^\ell a_im_i = 0,\quad \sum_{i=1}^\ell a_i = 0,\tag9.1
$$
consider the differential operator
$$
\square_{\bold{a}} = \prod_{a_j > 0} \left(\frac{\partial}{\partial c_j}
\right)^{a_j} - \prod_{a_j < 0} \left(\frac{\partial}{\partial c_j}
\right)^{-a_j}\tag9.2
$$
where we regard $c_1,\dots,c_\ell$ as variables on $L(\Cal{A}) = {\Bbb
C}^\ell$.  We will also consider the differential operators
$$
\Cal{Z}_0 = \sum_{j=1}^\ell c_j
\frac{\partial}{\partial c_j}, \quad
	\Cal{Z}_i = \sum_{j=1}^\ell m_{ji} c_j
\frac{\partial}{\partial c_j},\quad i = 1,\dots,n,
$$
where $m_j = (m_{j1},\dots,m_{jn}) \in {\Bbb Z}^n$ are the elements of
$\Cal{A}$.  Then, for a function $\Phi(c_1,\dots,c_\ell)$ on ${\Bbb
C}^\ell$, the {\it $\Cal{A}$-hypergeometric system with exponents
$\beta_0,\dots,\beta_n \in {\Bbb C}$\/} is the system of differential
equations
$$
\aligned \square_{\bold{a}} \Phi &= 0,\quad \hbox{for all $\bold{a}$
satisfying \thetag{9.1}}\\
	\Cal{Z}_i\Phi &= \beta_i \Phi,\quad \hbox{for all $i =
0,\dots,n$.} 
\endaligned\tag9.3
$$

	This system is holonomic, so that its solutions form a locally
constant sheaf outside a hypersurface in $L(\Cal{A}) = {\Bbb C}^\ell$.
Two especially nice facts are first, that the generic number of
linearly independent solutions of \thetag{9.3} is the normalized
volume of the polytope $\roman{Conv}(\Cal{A})$, and second, that
$(c_1,\dots,c_\ell)$ is generic if and only if the principal
$\Cal{A}$-determinant is nonvanishing, i.e., if and only if
$E_{\Cal{A}}(f) \ne 0$ for $f = \sum_{i=0} c_i \bold{t}^{m_i}$.
Proofs can be found in \cite{GKZ3}.

	A more direct connection with toric varieties can be seen by
considering the symbols of the operators $\square_{\bold{a}}$.  Using
$x_1,\dots,x_\ell$ as variables on the dual of $L(\Cal{A}) = {\Bbb
C}^\ell$, the symbol of $\square_{\bold{a}}$ is its Fourier transform
${\check \square}_{\bold{a}}$, which is obtained by replacing
$\partial/\partial c_j$ in \thetag{9.2} with $x_j$.  If we write
$\bold{a} = \bold{a}^+ - \bold{a}^-$ as in Lemma 5.5, then ${\check
\square}_{\bold{a}}$ becomes $x^{\bold{a}^+} - x^{\bold{a}^-}$.  Since
we do this for all $\bold{a}$ satisfying \thetag{9.1}, it follows from
\thetag{7.5} that we get the toric ideal $I_{\Cal{A}}$.  Hence the
Fourier transform of \thetag{9.3} is supported on the toric variety
$X_{\Cal{A}}$.  This toric variety plays an important role in the
proofs of many results about $\Cal{A}$-hypergeometric functions.
Also, although \thetag{9.3} involves infinitely many equations
$\square_{\bold{a}}\Phi =0$, one can always reduce to a finite number
using a Gr\"obner basis for the toric ideal $I_{\Cal{A}}$ (see
\cite{HLY1}).

	There are several ways to write down solutions to the
$\Cal{A}$-hypergeometric system.  For us, the most interesting method
involves the periods of the affine hypersurface $Z_f \subset ({\Bbb
C}^*)^n$ defined by $f = \sum_{i=0} c_i \bold{t}^{m_i}$.  We will
assume that $\Cal{A} = \Delta\cap{\Bbb Z}^n$, where $\Delta \subset
{\Bbb R}^n$ is an integer polytope.  As in \S5, we have the ring
$S_\Delta \subset {\Bbb C}[t_0,t_1,\dots,t_n]$.  In \cite{Bat4{\rm,
Thm.~14.2}}, the following result is proved.

\proclaim{Proposition 9.3} Let $\Cal{A} = \Delta\cap{\Bbb Z}^n$ and
$t_0^k\bold{t}^m \in S_\Delta$ (so that $m \in k\Delta\cap{\Bbb Z}^n$).
Then, for a $n$-cycle $\gamma \in H_n(({\Bbb C}^*)^n - Z_f)$ and $f =
\sum_{i=0}^\ell c_i \bold{t}^{m_i} \in L(\Cal{A})$ nondegenerate, the
period integral
$$
\Pi(c_1,\dots,c_\ell) = \Pi(f) = \int_\gamma \frac{\bold{t}^m}{f^k}\,
\frac{dt_1}{t_1}\wedge \cdots \wedge \frac{dt_n}{t_n} 
$$
satisfies the $\Cal{A}$-hypergeometric system \thetag{9.3} for
exponents $(\beta_0,\dots,\beta_n) = (-k,-m)$.
\endproclaim

	By Proposition 9.2, $Z_f$ is nondegenerate if and only if
$E_{\Cal{A}}(f) \ne 0$.  This helps explain why the singular set of
the $\Cal{A}$-hypergeometric system is defined by $E_{\Cal{A}} = 0$.
The integral in Proposition 9.3 is an example of an {\it Euler
integral\/}, which is the main object of study in \cite{GKZ2}.

	It is also possible to write down series solutions of the
$\Cal{A}$-hypergeometric system.  Formally, these solutions look like
$$
\Phi_\gamma(c_1,\dots,c_\ell) = \sum_{\bold{a}} \frac{\bold{c}^{\gamma
+ \bold{a}}}{\Pi_{j=1}^\ell \Gamma(\gamma_j + a_j + 1)},\tag9.4
$$
where $\gamma \in {\Bbb C}^\ell$, $\Gamma$ is the usual
$\Gamma$-function, and the sum is over all $\bold{a}$ satisfying
\thetag{9.1}.  By restricting to certain choices of $\gamma$ and 
$\beta_0,\dots,\beta_n$, one can show that the above series converge
locally when $(c_1,\dots,c_\ell)$ lies in certain regions of ${\Bbb
C}^\ell$ determined by a regular triangulation $T$ of
$\roman{Conv}(\Cal{A})$.  This is described as follows.  We saw in \S7
that a regular triangulation gives a maximal cone of $A_{\Cal{A}}$ in
the secondary fan of $\Cal{A}$.  Using the exact sequence
\thetag{7.3}, we get a cone $C(T) \subset \bigoplus_{m\in \Cal{A}}
{\Bbb R}\cdot e_m \simeq {\Bbb R}^\ell$.  Then one can prove that the
series \thetag{9.4} converges locally for those $(c_1,\dots,c_\ell)$
such that $-\log(c_1,\dots,c_\ell)$ (as defined in \thetag{6.1}) lies
in a suitable translate of $C(T)$.  Furthermore, one gets all
solutions of \thetag{9.3} in this way.  See \cite{GKZ3} for details.

	As mentioned earlier, the Newton polytope of the principal
$\Cal{A}$-determinant $E_{\Cal{A}}$ is the secondary polytope (with
vertices corresponding to regular triangulations).  In terms of
\thetag{6.1}, the regions of nice convergence of the series
\thetag{9.4} correspond (up to sign) to unbounded components of ${\Bbb
R}^\ell - \log(\{E_{\Cal{A}} = 0\})$.  This again illustrates why the
singular set of $\Cal{A}$-hypergeometric system is given by
$E_{\Cal{A}} = 0$.

	Finally, Proposition 9.3 shows that the system \thetag{9.3} is
closely tied to the torus.  Roughly speaking, the equations
$\Cal{Z}_i\Pi =\beta_i\Pi$ express the invariance of the period
integral $\Pi$ for $Z_f \subset ({\Bbb C}^*)^n$ under the
infinitesimal automorphisms of $({\Bbb C}^*)^n$.  If we were to
formulate a similar period integral for 
$Y_f \subset X_\Delta$, we would need to add further equations to
\thetag{9.3} to account for automorphisms of $X_\Delta$. 
This {\it extended $\Cal{A}$-hypergeometric system\/}
is described in \cite{HLY1}.

\head \S10. Mirror symmetry\endhead

	In 1991, a group of physicists made some startling enumerative
predictions concerning rational curves on a quintic threefold (see
\cite{CdGP}).  The basic idea was that a hard computation on the
quintic threefold became easier by working on its ``mirror''.  The
first mirror constructions involved finite quotients of weighted
projective spaces which, as Roan pointed out in 1992, can be described
naturally using toric methods (see \cite{Roa2}).  After Batyrev's 1993
introduction of reflexive polytopes (see \cite{Bat1}), toric geometry
has become a basic tool of mirror symmetry.

	We will not discuss mirror symmetry in general---the reader
should consult \cite{Yau}, \cite{Mor3} or \cite{CK} for an
introduction to this fascinating topic.  Rather, we will concentrate
on describing how toric geometry is used in mirror symmetry.  This is
a very active field, and the references given below are far from
complete.

	\subhead Complex and K\"ahler moduli\endsubhead The complex
moduli of varieties have been studied for many years.  The
infinitesimal deformations of a compact complex manifold $Y$ of
dimension $d$ are given by $H^1(Y,\Theta_Y)$, which in nice cases is
the tangent space to the complex moduli space.  When $Y$ is
Calabi-Yau, the isomorphism $\Omega_Y^d \simeq \Cal{O}_Y$ implies
$\Theta_Y \simeq \Omega_Y^{d-1}$, so that $H^1(Y,\Theta_Y)
\simeq H^{d-1,1}(Y)$.

	The intense study of K\"ahler moduli is more recent.  The
basic definition is as follows: given a K\"ahler manifold $Y$, the
{\it complexified K\"ahler cone\/} is the cone
$$
K_{\Bbb C}(Y) = \{B + iJ \in H^2(Y,{\Bbb C}) : J\ \hbox{is K\"ahler}\}
/ \roman{Im}\, H^2(Y,{\Bbb Z}),\tag10.1
$$
and the {\it complexified K\"ahler moduli space\/} is the quotient
$K_{\Bbb C}(Y)/\roman{Aut}(Y)$.  The exact structure of this quotient
space depends on how $\roman{Aut}(Y)$ acts on the ordinary K\"ahler
cone of $Y$ (as described in \S3).  In the Calabi-Yau case, it is
conjectured that the K\"ahler cone is polyhedral modulo the action of
$\roman{Aut}(Y)$ (see \cite{Mor1} for a precise statement).  This
would imply the existence of a semi-toric compactification of the
complexified K\"ahler moduli space.  Assuming that $\roman{Aut}(Y)$
acts discretely, the tangent space to the K\"ahler moduli space (we
will usually drop the adjective ``complexified'') is $H^{1,1}(Y)$.

	So far, we've assumed that $Y$ is smooth, and in dimension 3,
this is sufficient.  Higher dimensional generalizations of mirror
symmetry lead naturally to singular varieties, and one can define both
complex and K\"ahler moduli when $Y$ is a minimal Calabi-Yau variety
as in \S8.

	\subhead The naive idea of mirror symmetry\endsubhead A
Calabi-Yau threefold $Y$ together with a complexified K\"ahler class
$\omega = B + iJ$ determine a $N = 2$ superconformal field theory
(SCFT), and mirror symmetry suggests that there should be another
Calabi-Yau $Y^\circ$ with class $\omega^\circ$ which in some sense
interchanges the complex and K\"ahler structures of $Y$ and $Y^\circ$ but
still gives the same $N=2$ SCFT.

	This interchange of complex and K\"ahler moduli in particular
gives isomorphisms of the corresponding moduli spaces (the ``mirror
map'') and hence induces isomorphisms
$$
H^{d-1,1}(Y) \simeq H^{1,1}(Y^\circ),\quad H^{1,1}(Y) \simeq
H^{d-1,1}(Y^\circ)\tag10.2 
$$
between their tangent spaces.  But mirror symmetry is much more than
these isomorphisms, for an isomorphism of $N=2$ SCFTs also implies
that certain trilinear functions on $H^{1,1}(Y)$ and
$H^{d-1,1}(Y^\circ)$ should agree after a change of variables given by
the mirror map.  These trilinear functions are called {\it $3$-point
functions\/} or {\it correlation functions\/} in the physics
literature and are related to enumerative geometry and quantum
cohomology (for $H^{1,1}(Y)$) and Hodge theory (for
$H^{d-1,1}(Y^\circ)$).  Unfortunately, these topics are beyond the
scope of this survey.

	The physical theories involved in mirror symmetry have yet to
be defined rigorously, but there are ``mathematical mirror symmetry
conjectures'' which capture the mathematically interesting
consequences.  This can be done in all dimensions, and versions of
these conjectures are in \cite{Mor1--3}, \cite{Giv2}, \cite{Kont} and
\cite{Ver}.  Mirror symmetry for holomorphically symplectic manifolds
is proved in \cite{Ver}, and for complete intersections in projective
space, a version of mirror symmetry involving certain ``virtual''
numbers of rational curves of degree $d$ is proved in \cite{Giv1}.
Other than these exceptions, mathematical mirror symmetry is still
largely conjectural.  Hence, in what follows, our term ``mirror''
really means ``conjectural mirror''.

	\subhead Mirror symmetry for toric hypersurfaces\endsubhead In
\cite{Bat1}, Batyrev used reflexive polytopes to construct mirrors for
Calabi-Yau toric hypersurfaces.  As in \S8, the process starts with an
$n$-dimensional reflexive polytope $\Delta \subset M_{\Bbb R}$, which
by definition determines the anticanonical divisor on $X_\Delta$.
This gives canonical Calabi-Yau hypersurfaces $Y \subset X_\Delta$.
Then, if $\Sigma$ is a projective simplicial fan refining the normal
fan of $\Delta$ and satisfying $\Sigma(1) = N\cap\Delta^\circ -
\{0\}$, we obtain a family of minimal Calabi-Yau hypersurfaces
$\widehat Y \subset X_\Sigma$.

	To construct the mirror of this family, we repeat the above
procedure using the dual polytope $\Delta^\circ$.  Thus, for a
suitable fan $\Sigma^\circ$ in $M_{\Bbb R}$ with $\Sigma^\circ(1) =
M\cap\Delta - \{0\}$, the anticanonical divisor determines a family
$\widehat Y^\circ \subset X_{\Sigma^\circ}$ of minimal Calabi-Yaus.
This is conjectured to be the mirror of $\widehat Y \subset X_\Sigma$.

	Evidence for the mirror relation between $\widehat Y$ and
$\widehat Y^\circ$ comes from Theorem 8.2, which shows that
$$
H^{n-2,1}(\widehat Y) \simeq H^{1,1}(\widehat Y^\circ),\quad
H^{1,1}(\widehat Y) \simeq H^{n-2,1}(\widehat Y^\circ),
$$
as predicted by \thetag{10.2} (since $\widehat Y$ and $\widehat
Y^\circ$ have dimension $n-1$).  Also, $H^{n-2,1}(\widehat Y)$ has a
subspace $H^{n-2,1}_{\roman{poly}}(\widehat Y)$ consisting of
deformations of $\widehat Y$ obtained by varying its defining equation
in $X_{\Sigma}$, and similarly, $H^{1,1}(\widehat Y)$ has a subspace
$H^{1,1}_{\roman{toric}}(\widehat Y)$ consisting of restrictions of
$(1,1)$-classes on $X_{\Sigma}$.  The {\it monomial-divisor mirror
map\/} of \cite{AGM2} gives natural isomorphisms
$$
H^{n-2,1}_{\roman{poly}}(\widehat Y) \simeq
H^{1,1}_{\roman{toric}}(\widehat Y^\circ),\quad
H^{1,1}_{\roman{toric}}(\widehat Y) \simeq
H^{n-2,1}_{\roman{poly}}(\widehat Y^\circ).
$$

\demo{Example} To compute the mirror of the quintic threefold $Y
\subset {\Bbb P}^4$, note that $Y$ is an anticanonical hypersurface,
and the corresponding polytope $\Delta \subset M_{\Bbb R} \simeq {\Bbb
Z}^4$ is reflexive.  Hence the mirror should be determined by the dual
polytope $\Delta^\circ \subset N_{\Bbb R}$.

	The toric variety $X_{\Delta^\circ}$ comes from the normal fan
of $\Delta^\circ$.  The cone generators of this fan are the vertices of
$\Delta$, which are
$$
(-1,-1,-1,-1),\, (4,-1,-1,-1),\, (-1,4,-1,-1),\,
(-1,-1,4,-1),\, (-1,-1,-1,4).
$$
These generate a sublattice $M' \subset M$ of index 125, and
the quotient $M/M'$ is 
$$
H = \big\{(a_0,a_1,a_2,a_3,a_4) \in ({\Bbb Z}_5)^5 :
{\textstyle{\sum_{i=0}^4}} a_i \equiv 0 \bmod 5\big\}/{\Bbb Z}_5
$$
where ${\Bbb Z}_5 \subset ({\Bbb Z}_5)^5$ is the diagonal subgroup.
Using the lattice $M'$, the normal fan gives ${\Bbb P}^4$, so by
\cite{Oda1{\rm, Cor.~1.16}}, $X_{\Delta^\circ}$ is the quotient of
${\Bbb P}^4$ by $M/M' \simeq H$, where $[a_0,a_1,a_2,a_3,a_4] \in H$
acts on ${\Bbb P}^4$ by the map $(x_0,x_1,x_2,x_3,x_4) \mapsto
(\zeta^{a_0} x_0,\zeta^{a_1} x_1, \zeta^{a_2} x_2,\zeta^{a_3}
x_3,\zeta^{a_4} x_4)$ for $\zeta = \exp(2\pi i/5)$.

	The homogeneous coordinate ring $S$ of $X_{\Delta^\circ}$ from
\S2 is the polynomial ring $S= {\Bbb C}[x_0,x_1,x_2,x_3,x_4]$, which is
graded by the Chow group $A_3(X_{\Delta^\circ})$.  From \thetag{2.2},
we get
$$
0 \longrightarrow {\Bbb Z}^4 \overset\alpha\to\longrightarrow {\Bbb
Z}^5 \overset\beta\to\longrightarrow {\Bbb Z}\oplus H \longrightarrow
0,
$$
where $\alpha$ is as usual and $\beta$ is given by
$$
\beta(a_0,a_1,a_2,a_3,a_4) = \big({\textstyle{\sum_{i=0}^4}} a_i,
[-a_1-a_2-a_3-a_4,a_1,a_2,a_3,a_4]\big) \in {\Bbb Z}\oplus H.
$$
Thus $A_3(X_{\Delta^\circ}) \simeq {\Bbb Z}\oplus H$ and the grading
on $S$ is obtained by letting a monomial $x_0^{a_0} x_1^{a_1}
x_2^{a_2} x_3^{a_3} x_4^{a_4}$ have
degree $\beta(a_0,a_1,a_2,a_3,a_4) \in {\Bbb Z}\oplus H$.

	It follows that the anticanonical class has degree
$\beta(1,1,1,1,1) = (5,0)$, and the only monomials in $S$ of this
degree are $x_i^5$ and $x_0x_1x_2x_3x_4$ (this can be seen directly or
using the isomorphism $S_{(5,0)} \simeq L(\Delta^\circ)$ from
\thetag{6.2}).  Furthermore, the automorphisms of $X_{\Delta^\circ}$
given by its torus show that any anticanonical hypersurface is
isomorphic to one defined by an equation of the form
$$ 
x_0^5 + x_1^5 + x_2^5 + x_3^5 + x_4^5 + 5 \psi\, x_0 x_1 x_2 x_3 x_4 =
0.
$$
This gives a $1$-dimensional family of hypersurfaces $Y^\circ \subset
X_{\Delta^\circ}$.  In concrete terms, $Y^\circ$ is the quotient by
$H$ of the hypersurface in ${\Bbb P}^4$ defined by the above equation.

	Finally, to get the mirror, we pick a projective simplicial
fan $\Sigma^\circ$ in $M_{\Bbb R}$ refining the normal fan of
$\Delta^\circ$ such that $\Sigma^\circ(1) = M\cap\Delta-\{0\}$.  Then
the mirror family of the quintic threefold $Y$ is given by the
hypersurfaces $\widehat Y^\circ \subset X_{\Sigma^\circ}$ which are
the proper transforms of $Y^\circ \subset X_{\Delta^\circ}$.  Since $Y
\subset {\Bbb P}^4$ has $1$-dimensional K\"ahler moduli (i.e.,
$h^{1,1} = 1$), we expect its mirror $\widehat Y^\circ$ to have
$1$-dimensional complex moduli.  Note that $M\cap\Delta-\{0\}$ has
{\it lots\/} of points besides the vertices, and hence there are {\it
many\/} choices for $\Sigma^\circ$.  A specific choice for
$\Sigma^\circ$ is described in the appendix to
\cite{Mor3}.
\enddemo

	In general, once the mirror family has been found, one gets
enumerative predictions on $\widehat Y$ by computing the {\it mirror
map\/} and the {\it Yukawa coupling\/} on $\widehat Y^\circ$.  Without
getting into precise definitions, these objects are closely related to
period integrals and Picard-Fuchs equations on a toric hypersurface,
and the computations involve series expansions about certain
``special'' boundary points in the complex moduli space.  These
``special'' points have {\it maximal unipotent monodromy\/} and their
existence is suggested by the structure of the K\"ahler moduli space
of the mirror.  See \cite{Mor2-3} for further details.

	In many cases, the period integrals and Picard-Fuchs equations
can be computed directly (see \cite{Mor4} for some examples), though
Proposition 9.3 suggests a connection with $\Cal{A}$-hypergeometric
equations.  This has been studied carefully in \cite{HLY1-2}.  Further
references for both of these methods can be found in \cite{Mor2}.

\subhead Multiple mirrors and global moduli\endsubhead The theory
described so far is only local.  Given a family of minimal Calabi-Yau
hypersurfaces $\widehat Y \subset X_\Sigma$ and its mirror $\widehat
Y^\circ \subset X_{\Sigma^\circ}$, the complex moduli of $\widehat Y$
should correspond to the K\"ahler moduli of $\widehat Y^\circ$.  The
problem is that complex moduli typically form a quasi-projective
variety, while K\"ahler moduli are often the quotient of a bounded
domain by a finite group.  These are very different types of
mathematical objects.

	Hence, the symmetry $\widehat Y \leftrightarrow \widehat
Y^\circ$ only gives a local isomorphism between complex and K\"ahler
moduli.  To get a global version of mirror symmetry for toric
hypersurfaces, we need to examine the role of the many fans $\Sigma$
and $\Sigma^\circ$ which can occur.  This is the {\it multiple
mirror\/} phenomenon.  For complex moduli, we expect the $\widehat
Y$'s for different $\Sigma$'s to be related by a series of flops, and
hence they should have the same complex moduli (although the varieties
themselves are not isomorphic).  

	The picture on the K\"ahler side is more interesting, for
here, we've seen in \S7 that the various K\"ahler cones for
$\Sigma^\circ(1) = M\cap\Delta-\{0\}$ fit together to form the
secondary fan of $\Cal{B} = M\cap\Delta-\{0\}$.  The support of this
fan is a strongly convex cone in $A_{n-1}(X_\Delta)\otimes{\Bbb R}$.
Gluing together the corresponding complexified K\"ahler cones
\thetag{10.1} gives the {\it partially enlarged K\"ahler moduli
space\/} from \cite{AGM1}.  (Strictly speaking, we are only dealing
with the toric part of the moduli space, but we will ignore this
detail.)

	The partially enlarged moduli space is still too small to be
quasi-projective.  To get something bigger, we use the enlarged
secondary fan from \S7.  This is determined by an ample divisor on a
toric variety, which here is the anticanonical divisor.  The
corresponding fan is the secondary fan $\Cal{N}(\Sigma(\Cal{A}))$ for
$\Cal{A} = M\cap\Delta$.  When we complexify the cones in this fan and
glue them together, we get the {\it enlarged K\"ahler moduli space\/}
from \cite{AGM1}, which corresponds to the whole complex moduli of the
mirror.

	This isomorphism between moduli spaces extends to certain
natural compactifications.  For simplicity, let's restrict to
polynomial moduli, which come from those $f \in L(\Delta)$ which are
nondegenerate (i.e., $E_{\Cal{A}}(f) = 0$ for $\Cal{A}$ as above)
modulo automorphisms of the toric variety.  If we instead use only
automorphisms coming from the torus, we can compactify using the
Newton polytope of the principal $\Cal{A}$-determinant $E_{\Cal{A}}$.
The resulting compactification has a natural toric structure where the
fixed points correspond to the vertices of the Newton polytope.  But
the Newton polytope of $E_{\Cal{A}}$ is the secondary polytope
$\Sigma(\Cal{A})$, whose vertices correspond to cones in the enlarged
K\"ahler moduli space.  Furthermore, the fixed points for vertices
corresponding to cones in the partially enlarged K\"ahler moduli space
are precisely the ``special'' points mentioned earlier.  For more of
the mathematics behind this picture, see \cite{AGM2}.

	The physics interpretation is also quite interesting: the
cones coming from the partially enlarged K\"ahler moduli space
correspond to different ``Calabi-Yau phases'' of the same SCFT, while
the other cones in the enlarged moduli space correspond to
``non-geometric phases''.  See \cite{AGM1} and \cite{MP} for more
details.

	\subhead Nef-partitions and Gorenstein cones\endsubhead In
\S8, we described how a nef-partition $\roman{Vert}(\Delta^\circ) =
E_1\cup \cdots \cup E_r$ of a reflexive polytope $\Delta$ gave a
Minkowski sum $\Delta = \Delta_1 + \cdots + \Delta_r$ in $M_{\Bbb R}$
and a canonical Calabi-Yau complete intersection $V =
Y_1\cap\cdots\cap Y_r \subset X_\Delta$.  Furthermore, a projective
simplicial fan $\Sigma$ refining the normal fan of $\Delta$ and
satisfying $\Sigma(1) = N\cap\Delta^\circ-\{0\}$ gives the
MPCP-desingularization $\widehat V = \widehat Y_1\cap\cdots\cap
\widehat Y_r \subset X_{\Sigma}$ as in \thetag{8.5}.  This is a family
of minimal Calabi-Yau complete intersections.

	To get the mirror family, we follow \cite{LBor} and consider
the polytopes
$$
\nabla_i = \roman{Conv}(\{0\}\cup E_i) \subset N_{\Bbb R},\quad i =
1,\dots,r.\tag10.3 
$$
Then $\nabla = \nabla_1 + \cdots + \nabla_r$ is a reflexive polytope
in $N_{\Bbb R}$ with a natural nef-partition coming from
\thetag{10.3}.  It is interesting to note that $\nabla^\circ \subset
M_{\Bbb R}$ is different from our original $\Delta$.  In fact,
\cite{LBor} shows that
$$
\nabla^\circ = \roman{Conv}(\Delta_1,\dots,\Delta_r) \subset \Delta_1
+ \cdots + \Delta_r = \Delta \subset M_{\Bbb R},
$$
with a similar relation between $\Delta^\circ$ and $\nabla$ in
$N_{\Bbb R}$. 

	Then the nef-partition for $\nabla$ gives a family of
canonical Calabi-Yau complete intersections $V^\circ \subset
X_{\nabla}$, and picking a projective simplicial fan $\Sigma^\circ$
refining the normal fan of $\nabla$ and satisfying $\Sigma^\circ(1) =
M\cap\nabla^\circ-\{0\}$ gives a family $\widehat V^\circ \subset
X_{\Sigma^\circ}$ of minimal Calabi-Yau complete intersections which
is conjectured in \cite{LBor} to be the mirror of $\widehat V \subset
X_\Delta$.  Evidence for this is presented in \cite{BB2-3},
\cite{Bv} and \cite{LT} (see also the references in \cite{BB3}).

	A significant generalization of this construction, which uses
{\it Gorenstein cones\/}, appeared in \cite{BB1}.  We will not go into
the details, but Gorenstein cones can explain all of the mirror
constructions we've given so far, as well as describing mirrors for
certain rigid Calabi-Yau manifolds (where the mirror may have a
different dimension).  Details can be found in
\cite{BB1}. 

	\subhead Reid's fantasy\endsubhead If a Calabi-Yau threefold
$Y_1$ degenerates to a variety with only nodes as singularities, then
one can resolve the singularities to obtain another Calabi-Yau $Y_2$
(which could fail to be K\"ahler).  Locally, a vanishing cycle $S^3
\subset Y_1$ collapses to a node and is then resolved to give ${\Bbb
P}^1 \simeq S^2 \subset Y_2$.  Hence $Y_1$ and $Y_2$ can have quite
different Betti numbers.  In \cite{Rei3}, Reid speculates that the
moduli of {\it all\/} Calabi-Yau threefolds many be connected in this
(possibly non-K\"ahler) way.

	In the physics literature, the singular Calabi-Yau between
$Y_1$ and $Y_2$ is a {\it conifold\/} and going from $Y_1$ to $Y_2$ is
a {\it conifold transition\/}.  Such transitions were long thought to
produce unacceptable singularities in the physical theories, but
recently (see \cite{GMS}), these difficulties were resolved by
allowing certain non-perturbative string states (electrically charged
black holes) on $Y_1$ to become massless on the conifold and to be
intrepreted on $Y_2$ as elementary perturbative states (elementary
particles).  It follows that a K\"ahler version of Reid's fantasy
would enable any two Calabi-Yau threefolds to be connected by a single
physical theory.  This implies that when ${\Bbb R}^4 \times
(\hbox{Calabi-Yau threefold})$ is used to model the vacuum state of
the universe, we don't need to worry about which Calabi-Yau to use,
since all can occur.

	Although Reid's fantasy is still conjectural, it has been
verified that for the $7555$ Calabi-Yau threefolds mentioned in \S8,
their moduli are connected through K\"ahler varieties, though the
singularities may be more complicated than nodes (see \cite{ACJM} and
\cite{CGGK}).  The basic idea is as follows.  Suppose we have
reflexive polytopes $\Delta_2 \subset \Delta_1$, which implies
$\roman{Vert}(\Delta_2) \subset \roman{Vert}(\Delta_1)$.  Then, for
$Y_1 \subset X_{\Delta_1}$ defined by $f \in L(\Delta_1)$, we can
degenerate $Y_1$ by letting the coefficients of $f$ corresponding to
vertices not in $\Delta_2$ become zero.  This gives a singular variety
$\widetilde Y \subset X_{\Delta_1}$ which, when resolved, corresponds
to a Calabi-Yau hypersurface $Y_2 \subset X_{\Delta_2}$.  Once the
MPCP-desingularizations are taken into account, we can link the
corresponding moduli spaces.  As already mentioned, the singularities
of $\widetilde Y$ may be more complicated than just simple nodes, and
at present there is no physical explanation of the transition from
$Y_1$ to $Y_2$.  But this method is sufficient to link up the $7555$
Calabi-Yaus on the list, and as noted in \cite{Ska}, it may be
sufficient to connect the moduli of all $3$-dimensional Calabi-Yau
toric hypersurfaces.

	\subhead Further remarks\endsubhead Although our discussion of
mirror symmetry has been rather superficial (and has omitted some
important ideas), it should be clear that toric geometry has a
prominent role to play, if for no other reason than providing a rich
supply of examples.  Notice also that virtually everything in the
earlier sections of the paper has been used.  Symplectic geometry
seems to be an exception, but this is only because we didn't describe
Witten's {\it linear sigma models\/} \cite{Wit}, which are physical
theories where toric varieties enter by means of symplectic reduction.
See \cite{MP} for more details and for some other interesting uses of
toric geometry in mirror symmetry.

	We should also mention that it is possible to compute the
quantum cohomology of toric varieties (see \cite{Bat3}).  In addition,
certain mirror symmetry calculations suggest that in some situations,
the usual Hodge numbers $h^{p,q}$ need to be replaced by {\it
string-theoretic Hodge numbers\/} $h_{\roman{st}}^{p,q}$.  For
Calabi-Yau complete intersections in a toric varieties, these numbers
are computed in \cite{BB2}.

\head \S11. Other developments\endhead

	Besides the topics reported on so far, the last few years have
seen a lot of interesting work on other aspects of toric geometry.
Here is a selection, with apologies for the many fine papers not
mentioned.

	\subhead Very ample divisors\endsubhead As is well-known, an
ample divisor D on a complete toric variety $X$ is very ample if $X$
is smooth or has dimension $\le 2$.  In general, $D$ may fail to be
very ample, but \cite{EW} proves that $(n-1)D$ is always very ample,
where $n \ge 2$ is the dimension of $X$.  Examples show that this
result is sharp.

	There is also a notion of {\it $k$-very ampleness\/} which
measures the behavior of $D$ relative to $0$-dimensional subschemes $Z
\subset X$ with $h^0(\Cal{O}_Z) = k+1$ (so that $0$-ample means
spanned by global sections and $1$-ample means very ample).  For
smooth toric surfaces, \cite{DiR} shows that $k$-very ampleness can be
interpreted in terms of convexity properties of the support function
of $D$ relative to the integer lattice.

	\subhead Embeddings into toric varieties\endsubhead Another
nice fact about toric varieties concerns embeddings of a complete
variety $Y$.  The Chevalley criterion (proved by Kleiman) states that
$Y$ can be embedded into a projective space if and only if every
finite subset of $Y$ is contained in an affine open.  When $Y$ is
normal, \cite{W\l o} proves that we can embed $Y$ into a complete
toric variety if and only if every {\it two element\/} subset of $Y$
lies in an affine open.

	\subhead Classifying toric varieties\endsubhead There are
several ways one can try to classify toric varieties.  For example,
one can work one dimension at a time, which is the approach taken in
\S8 when we discussed reflexive polytopes.  Another strategy is to
classify smooth complete toric varieties according to their Picard
number $\rho$.  For $\rho = 2$, this is done in \cite{Kle}, which such
varieties are shown to be projective (and \cite{ESch} finds projective
embeddings for which Conjecture 5.6 is satisfied.)  Similarly,
\cite{KS} shows that smooth toric varieties with Picard number $3$ are
projective.  Using this and the primitive collections defined in \S2,
a classification for smooth complete toric varieties with $\rho = 3$ is
given in \cite{Bat2}.  The papers mentioned here also contain
references to earlier work on classification.

	\subhead Invariants of toric varieties\endsubhead When a toric
variety $X$ is smooth (resp.\ simplicial), its cohomology over ${\Bbb
Z}$ (resp.\ ${\Bbb Q}$) is well-known.  However, it is also possible
to compute the rational intersection cohomology of a compact toric
variety (see \cite{Fie}), and since a toric variety has a natural
torus action, there is also equivariant cohomology to consider, which
is computed in \cite{Bif}.

	Turning to less topological invariants, one can study the
K-theory (Grothendieck groups of vector bundles or coherent sheaves)
of a toric variety, which coincide in this case and are computed in
\cite{More1}.  The Brauer group of a toric variety is discussed in
\cite{DFM} in the case of an algebraically closed field, and the split
case over arbitrary fields is studied in \cite{For}.  It would be
interesting to see what results could be obtained for non-split tori
(to be defined below).  The paper \cite{For} also considers certain
invariants of $X_\Sigma$ which depend only on the combinatorial type
of the fan $\Sigma$.  Similar questions about $\roman{Pic}(X_\Sigma)$
are studied in \cite{Eik}.

	Another invariant known in the smooth (resp.\ simplicial)
toric case is the Chow ring $A^*(X)$ (resp.\ $A^*(X)\otimes{\Bbb Q}$).
For a general toric variety, the Chow groups $A_k(X)$ are computed in
\cite{Dan2{\rm, Sect.~10}}, and when $X$ is a complete toric variety,
one can relate $A_k(X)\otimes{\Bbb Q}$ to a certain weight filtration
of the Borel-Moore homology of $X$ (see \cite{Tot}).  Besides these
classical Chow groups, there are also the {\it operational Chow groups
$A^k(X)$\/} of Fulton and MacPherson, which give a Chow cohomology
ring $A^*(X)$.  When $X$ is a complete toric variety, this ring is
computed in \cite{FS}.

	\subhead Intersection theory on toric varieties\endsubhead As
mentioned in the introduction, Hard Lefschetz for simplicial toric
varieties was used in Stanley's proof of McMullen's conjectures about
convex simplicial polytopes.  A nice discussion of this may be found
in \cite{Ful{\rm, Sect.~5.6}}.  Since Hard Lefschetz for intersection
cohomology is a very deep result, it is reasonable to ask if a simpler
proof exists in the toric case.  This led to the papers \cite{Oda1}
and \cite{Oda4} on the de Rham cohomology of toric varieties, although
the question of finding a toric proof of Hard Lefschetz is still open.
(We should also mention the paper \cite{Dan1}, which studies the de
Rham cohomology of toriodal varieties.)

	Subsequently, a proof McMullen's conjectures which avoided
toric varieties and Hard Lefschetz was found by McMullen \cite{McM}.
His proof used a certain {\it polytope algebra\/}.  In \cite{FS}, it
is shown that in dimension $n$, the polytope algebra is the inverse
limit of the Chow cohomology rings $A^*(X)\otimes{\Bbb Q}$ over the
directed system of {\it all\/} toric compactifications of the torus
$({\Bbb C}^*)^n$.

	\subhead Counting lattice points\endsubhead If $X_\Delta$ is
the toric variety determined by an $n$-dimen\-sion\-al integer polytope
$\Delta \subset M_{\Bbb R} \simeq {\Bbb R}^n$, then its Todd class can
be written
$$
Td(X_\Delta) = \sum_{\sigma\in\Sigma} r_\sigma [V(\sigma)] \in
{\textstyle{\bigoplus}_{k=0}^n} A_k(X_\Delta),
$$
where $\Sigma$ is the normal fan of $\Delta$, $[V(\sigma)]$ is the
class of the orbit closure corresponding to $\sigma$, and $r_\sigma
\in {\Bbb Q}$.  By \cite{Dan2}, the number of integer points in
$\Delta$ is given by
$$
l(\Delta) = \sum_{\sigma\in\Sigma} r_\sigma
\roman{Vol}(F(\sigma)),
$$
where $\roman{Vol}(F(\sigma))$ is the normalized volume of the face of
$\Delta$ corresponding to $\sigma$.  There is a similar formula for
the Ehrhart polynomial $E_\Delta(k) = l(k\Delta)$ for $k \ge 1$.  See
\cite{Ful{\rm, Sect.~5.3}} or \cite{More2{\rm, Sect.~1.1}} for a nice
introduction to this topic.

	These formulas reduce the problem of counting lattice points
to finding an explicit expression for the Todd class.  While this can
be done for any given toric variety (see \cite{Ful}), general formulas
weren't available until recently, when several different solutions
were found.  In \cite{More2}, a formula for $r_\sigma$ is given which
depends {\it only\/} on the cone $\sigma$ and not the fan in which it
sits.  For simplicial toric varieties, \cite{Pom} introduces the idea
of a ``mock'' Todd class, denoted $TD(X_\Delta)$, which is built using
formulas from the smooth case.  The difference $Td(X_\Delta) -
TD(X_\Delta)$ is then described using functions which involve Dedekind
sums and lead to explicit formulas for the component $Td^2(X_\Delta)$.
(For the weighted projective space ${\Bbb P}(q_0,q_1,q_2,q_3)$, another
approach to this computation can be found in \cite{Lat}.)  A slighlty
different definition of ``mock'' Todd class is given in \cite{CS} and
is used to obtain a third expression for $Td(X_\Delta)$.  A corollary
is an explicit formula for the number of lattice points in an integer
simplex.  

	We should also mention that the Todd class is related to the
total Chern class.  For a singular toric variety, the formula
$ch(X_\Delta) = \sum_{\sigma\in\Sigma} [V(\sigma)]$ is shown to hold
in homology with closed supports (see \cite{BBF}).

	The quite different approach to the study of lattice points
appears in \cite{Bri1}.  Here, the object of interest is the Laurent
polynomial
$$
l_\Delta(\bold{t}) = \sum_{m\in \Delta\cap M}
\bold{t}^m, 
$$
where we still assume that $X_\Delta$ is simplicial.  Using the
Lefschetz-Riemann-Roch theorem from equivariant K-theory,
$l_\Delta(\bold{t})$ can be written as a sum of rational functions in
$\bold{t}$ determined by the cones at the vertices of $\Delta$.  Hence
we get a formula which not only counts lattice points (by setting
$\bold{t} = (1,\dots,1)$) but also describes the lattice points
themselves.  Simpler proofs of this result can be found in \cite{Ish}
and \cite{SI}, and a weighted version is in \cite{Bri2}.  We should
also mention that these questions can be studied from a purely
``polytope'' point of view, which uses a {\it combinatorial
Riemann-Roch theorem} and avoids toric methods.  See \cite{KK} for
details.

	\subhead Rational points on toric varieties over number
fields\endsubhead To define a toric variety over a number field $K$,
first observe that a torus $T$ over $K$ is determined by a lattice $M
\simeq {\Bbb Z}^n$ with an action of $\roman{Gal}(E/K)$, where $K
\subset E$ is a finite Galois extension over which the torus splits,
i.e., becomes isomorphic to $({\Bbb G}_m)^n$.  Then a toric variety
$X$ over $K$ containing $T$ is determined by a fan in $M_{\Bbb R}$
which is invariant under $\roman{Gal}(E/K)$.  See \cite{BT} for
details and references.

	In this situation, one can study $K$-rational points on $X$
using the height function coming from a metrized ample line bundle
$\Cal{L}$ on $X$.  The basic question concerns the asymptotics of
$N(T,\Cal{L},B)$, which is the number of $K$-rational points in the
torus $T \subset X$ with height bounded by $B$.  When $X$ is a smooth
projective Fano toric variety over $K$, the natural line bundle to use
is the anticanonical line bundle $\Cal{L} = \Cal{O}(-K_X)$.  In
\cite{BT}, the asymptotic formula
$$
N(T,\Cal{L},B) = c\,B\,(\log B)^{r-1}(1 + o(1)),\quad B \gg 0 
$$
is proved, where $r$ is the rank of the Picard group $\roman{Pic}(X)$
over $K$ and $c$ is a nonzero constant.  This verifies the toric case
of a conjecture of Manin for Fano varieties over number fields (see
\cite{FMT}).  The constant $c$ can be explicitly computed in terms of
the geometry of the cone of effective divisors in
$\roman{Pic}(X)_{\Bbb R}$, the order of the non-trivial part of the
Brauer group of $X$, and a certain Tamagawa number associated with the
metrized bundle $\Cal{L}$ on $X$.

	\subhead Quotients of toric varieties\endsubhead In \S2, we
constructed a toric variety $X$ as the quotient $({\Bbb C}^{\Delta(1)}
- Z(\Sigma))/G$, which is a geometric quotient when $X$ is simplicial.
Oda observed that ${\Bbb C}^{\Delta(1)} - Z(\Sigma)$ has a natural
structure of a toric variety, so that we are taking the quotient of a
toric variety by a subgroup of its torus.  One can study this problem
in general, and various types of quotients are possible, including
{\it combinatorial quotients\/}, {\it GIT quotients\/} and {\it Chow
quotients\/}.  For GIT quotients, there can be several quotients in
any given situation, and the different quotients are related by a
chamber structure where quotients for chambers sharing a common wall
are (in good cases) related by a flip.  (This is similar to the
secondary fan described in \S7.)  Furthermore, the Chow quotient can
be described using a certain fiber polytope and is the inverse limit
of the GIT quotients.  See \cite{B-BS}, \cite{Hu},
\cite{KSZ2} and \cite{Tha} for details.

	\subhead Residues on toric varieties\endsubhead There is a
huge literature on residues.  For multidimensional residues, toric
methods can be used in defining local residues (see \cite{PT}), and
Gr\"obner methods, including the Gr\"obner fan from \S7, can be used
in computing global residues (see \cite{CDS}).  A version of global
residues specific to toric varieties was introduced in \cite{Cox3},
and various properties of these {\it toric residues\/} have been
studied in \cite{CCD} and \cite{CD}.  Toric residues have also been
used in mirror symmetry (see \cite{MP}).

	\subhead Singularities of toric varieties\endsubhead
Singularities of affine toric varieties have been the subject of
several recent papers.  For example, basic tools in deformation theory
are the spaces $T^1_X$ and $T^2_X$ which describe infinitesimal
deformations and obstructions, and for an affine toric variety $X$,
these are computed in \cite{Alt1}.  One can also study deformations
which are themselves toric, and in the case of isolated
$3$-dimensional toric Gorenstein singularities, the versal deformation
can be constructed entirely by toric means (see \cite{Alt3}).
Furthermore, the irreducible components of the deformation space
correspond to certain decompositions of a lattice polytope into
Minkowski sums of other lattice polytopes.  The role of Minkowski sums
in the deformation theory of affine toric varieties is also explored
in \cite{Alt2}.  Another reference for Gorenstein toric singularities
is \cite{Nom}.

	In studying the resolutions of singularities of an affine
toric variety $X$, a divisor in a resolution of $X$ is {\it
essential\/} if a birational copy of the divisor appears in every
resolution of the variety, and it is {\it equivariant essential\/} if
it appears (again up to birational equivalence) in every toric
resolution.  In \cite{BG-S}, it is shown that if $X$ comes from a cone
$\sigma \subset N_{\Bbb R}$, then equivariant essential divisors
correspond bijectively to minimal generators of the semigroup
$\sigma\cap N - \{0\}$.  Furthermore, in dimension $3$, the same is
true for essential divisors.  Other papers dealing with the resolution
of $3$-dimensional toric singularities are \cite{ESpa} and \cite{Pou}.

	In another context, a conjecture of Shokurov on the minimal
discrepancies of log-terminal singularities was verified for toric
singularities in \cite{ABor}.  Hypersurface sections of toric
singularities are considered in \cite{Tsu}.  It is also possible to
discuss toric singularities without reference to a base field (or even
a base scheme).  This topic is studied in \cite{Kat} and may have
applications to arithmetic algebraic geometry.

	\subhead Resolution of singularities\endsubhead A recent
development is the discovery in \cite{Ad} and \cite{BP} of a simple
proof of a weak form of resolution of singularities.  The precise
result is that if $X$ is a normal projective variety in characteristic
$0$ and $Y\subset X$ is a proper subvariety, then there is a
birational morphism $f : \widehat X \to X$ such that $\widehat X$ is a
smooth projective variety and $f^{-1}(Y)$ is a strict normal-crossings
divisor.  However, $f : \widehat X - f^{-1}(Y) \to X - Y$ might not be
the identity, so we don't get a resolution of singularities in the
usual sense.  The proofs in \cite{Ad} and \cite{BP} are slightly
different, but both make essential use of toric methods.

\head Conclusion\endhead

	In this survey, we have attempted to convey the richness of
the recent work in toric geometry.  An unexpected consequence of all
this activity is that it is less clear where to learn about the subject.
One could start with the standard approach to toric varieties, as in
\cite{Dan2}, \cite{Ful}, \cite{Oda2}.  (Other introductions include
\cite{Ewa}, \cite{Oda3} and \cite{Rei4}.)  Alternatively, one could
begin with the quotient construction of \S2, where \cite{Cox2} is one
of many references.  This is closely tied to the symplectic approach,
as described in \cite{Aud} and \cite{Gui}.  Yet another starting point
would be the theory of non-normal toric varieties, as in \cite{GKZ1}
or \cite{Stu2}.

	Of course, there is no ``best'' approach to toric varieties.
The multiplicity of entry points is actually a virtue, for it enables
people from different areas of mathematics to learn about and
contribute to this fascinating and accessible part of algebraic
geometry.

	The author is grateful to Victor Batyrev, Sheldon Katz, David
Morrison and Bernd Sturmfels for numerous useful comments.  The
referee made many thoughtful suggestions which are greatly
appreciated.

\refstyle{A}

\enddocument